\documentclass[preprint]{elsarticle}

\usepackage{lineno,hyperref}
\usepackage{amssymb}
\usepackage{amsmath}
\usepackage{amsthm}
\theoremstyle{remark}
\newtheorem{remark}{Remark}

\theoremstyle{plain}
\newtheorem{proposition}{Proposition}

\theoremstyle{plain}
\newtheorem{lemma}{Lemma}

\theoremstyle{plain}
\newtheorem{theorem}{Theorem}

\theoremstyle{plain}
\newtheorem{definition}{Definition}

\theoremstyle{plain}
\newtheorem{question}{Question}

\theoremstyle{plain}
\newtheorem{corollary}{Corollary}

\newcommand{\Real}{\mathbb R}
\newcommand{\set}[1]{\left\{#1\right\}}
\newcommand{\abs}[1]{\left\vert#1\right\vert}
\newcommand{\real}[1]{{\mathbb R}^{#1}}
\newcommand{\norm}[1]{\left\Vert#1\right\Vert}

\newcommand{\bff}{{\boldsymbol f}}
\newcommand{\bg}{{\boldsymbol g}}

\newcommand{\bu}{{\boldsymbol u}}

\newcommand{\bw}{{\boldsymbol w}}
\newcommand{\bx}{\boldsymbol x}
\newcommand{\by}{{\boldsymbol y}}

\newcommand{\Lspace}{\mathbb{L}}
\newcommand{\bl}{{\boldsymbol l}}
\newcommand{\X}{\mathbb{X}}

\newcommand{\bzero}{{\bf 0}}

\modulolinenumbers[200]

\journal{Journal of \LaTeX\ Templates}
\begin{document}

\begin{frontmatter}

\title{An Optimal Control Theory for the Traveling Salesman Problem and Its Variants}

%\tnotetext[mytitlenote]{}

\author{I. M. Ross\fnref{rossftnt}}
\author{R. J. Proulx\fnref{proulxftnt}}
\author{M. Karpenko\fnref{karpftnt}}
\address{Naval Postgraduate School, Monterey, CA 93943}
\fntext[rossftnt]{Distinguished Professor \& Program Director, Control and Optimization }
\fntext[proulxftnt]{Research Professor, Control and Optimization Laboratories}
\fntext[karpftnt]{Research Associate Professor \& Associate Director, Control and Optimization Laboratories}

%%%% or include affiliations in footnotes:
%%\author[mymainaddress,mysecondaryaddress]{Elsevier Inc}
%%\ead[url]{www.elsevier.com}
%%
%%\author[mysecondaryaddress]{Global Customer Service\corref{mycorrespondingauthor}}
%%\cortext[mycorrespondingauthor]{Corresponding author}
%%\ead{support@elsevier.com}
%%
%%\address[mymainaddress]{1600 John F Kennedy Boulevard, Philadelphia}
%%\address[mysecondaryaddress]{360 Park Avenue South, New York}

\begin{abstract}
We show that the traveling salesman problem (TSP) and its many variants may be modeled as functional optimization problems over a graph.  In this formulation, all vertices and arcs of the graph are functionals; i.e., a mapping from a space of measurable functions to the field of real numbers. Many variants of the TSP, such as those with neighborhoods, with forbidden neighborhoods, with time-windows and with profits, can all be framed under this construct. In sharp contrast to their discrete-optimization counterparts, the modeling constructs presented in this paper represent a fundamentally new domain of analysis and computation for TSPs and their variants. Beyond its apparent mathematical unification of a class of problems in graph theory, the main advantage of the new approach is that it facilitates the modeling of certain application-specific problems in their home space of measurable functions.  Consequently, certain elements of economic system theory such as dynamical models and continuous-time cost/profit functionals can be directly incorporated in the new optimization problem formulation. Furthermore, subtour elimination constraints, prevalent in discrete optimization formulations, are naturally enforced through continuity requirements. The price for the new modeling framework is nonsmooth functionals. Although a number of theoretical issues remain open in the proposed mathematical framework, we demonstrate the computational viability of the new modeling constructs over a sample set of problems to illustrate the rapid production of end-to-end TSP solutions to extensively-constrained practical problems.

\end{abstract}


\end{frontmatter}

\linenumbers

\section{Introduction}
\begin{figure}[h!]
      \begin{center}
      \framebox{\parbox{3.0in}{
      \centering{\includegraphics[trim={1mm, 1mm, 1mm, 1mm}, clip, angle=0, width = 0.45\textwidth]{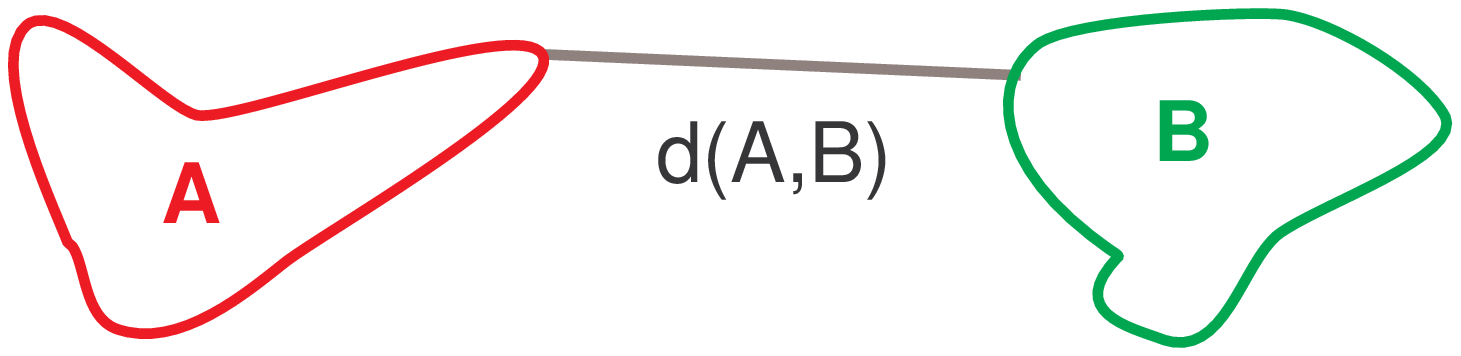}}
      }} %\vskip 5mm
      \caption{\textsf{A candidate distance function $d(A,B)$ between two disjoint sets, $A$ and $B$.}}
    \label{fig:2-City-Blobs}
    \end{center}
   \end{figure}

\begin{figure}[h!]
      \centering
      \framebox{\parbox{3.0in}{\centering
      {\includegraphics[trim={1mm, 1mm, 1mm, 1mm}, clip, angle=0, width = 0.45\textwidth]{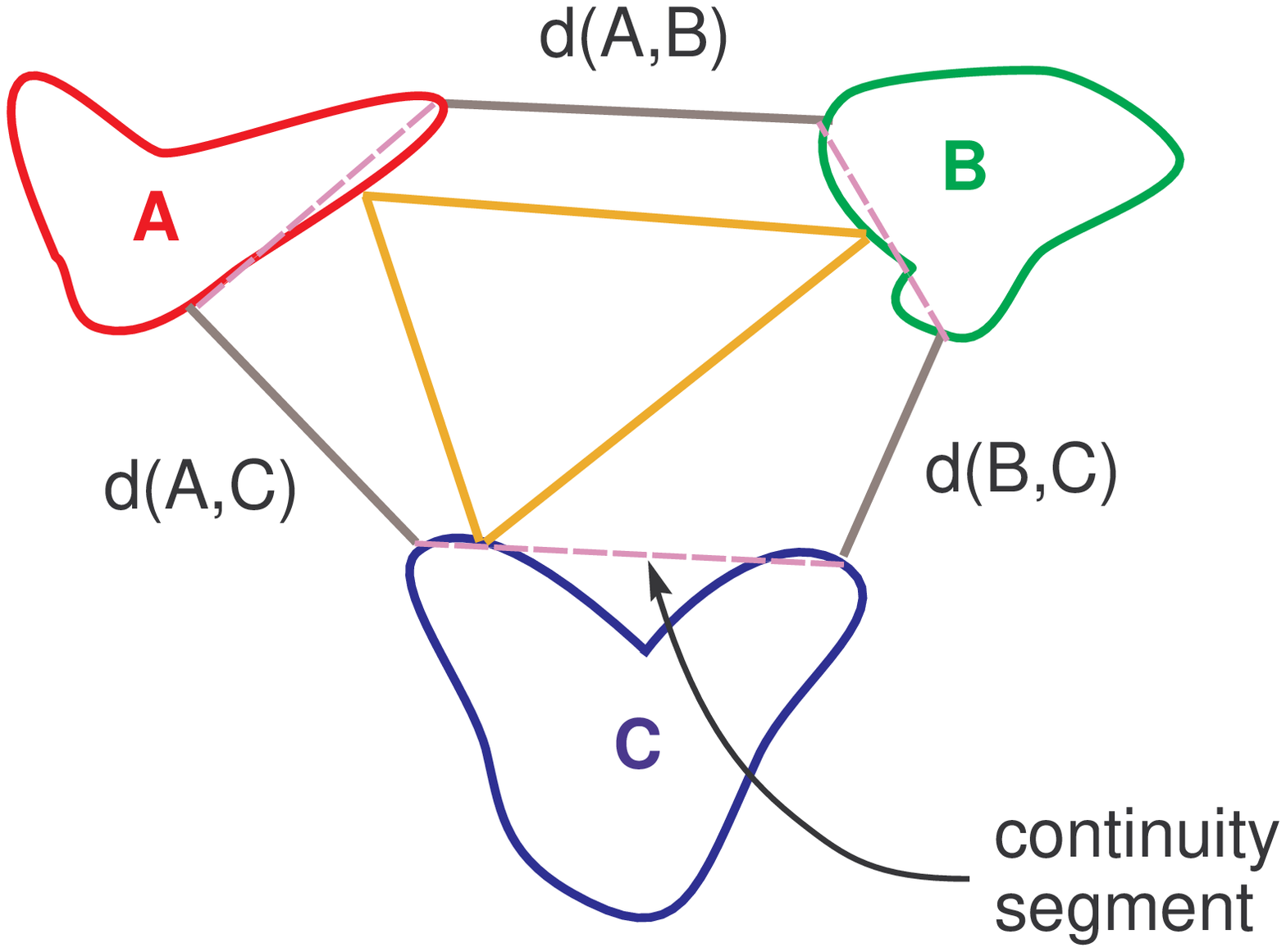}}
      }}
      \caption{\textsf{Illustrating how the shortest distance between two cities does not produce the shortest tour.}}
    \label{fig:3-City-Blobs}
   \end{figure}

\begin{figure}[h!]
      \centering
      \framebox{\parbox{3.0in}{\centering
      {\includegraphics[trim={1mm, 1mm, 1mm, 1mm}, clip, angle=-0, width = 0.5\textwidth]{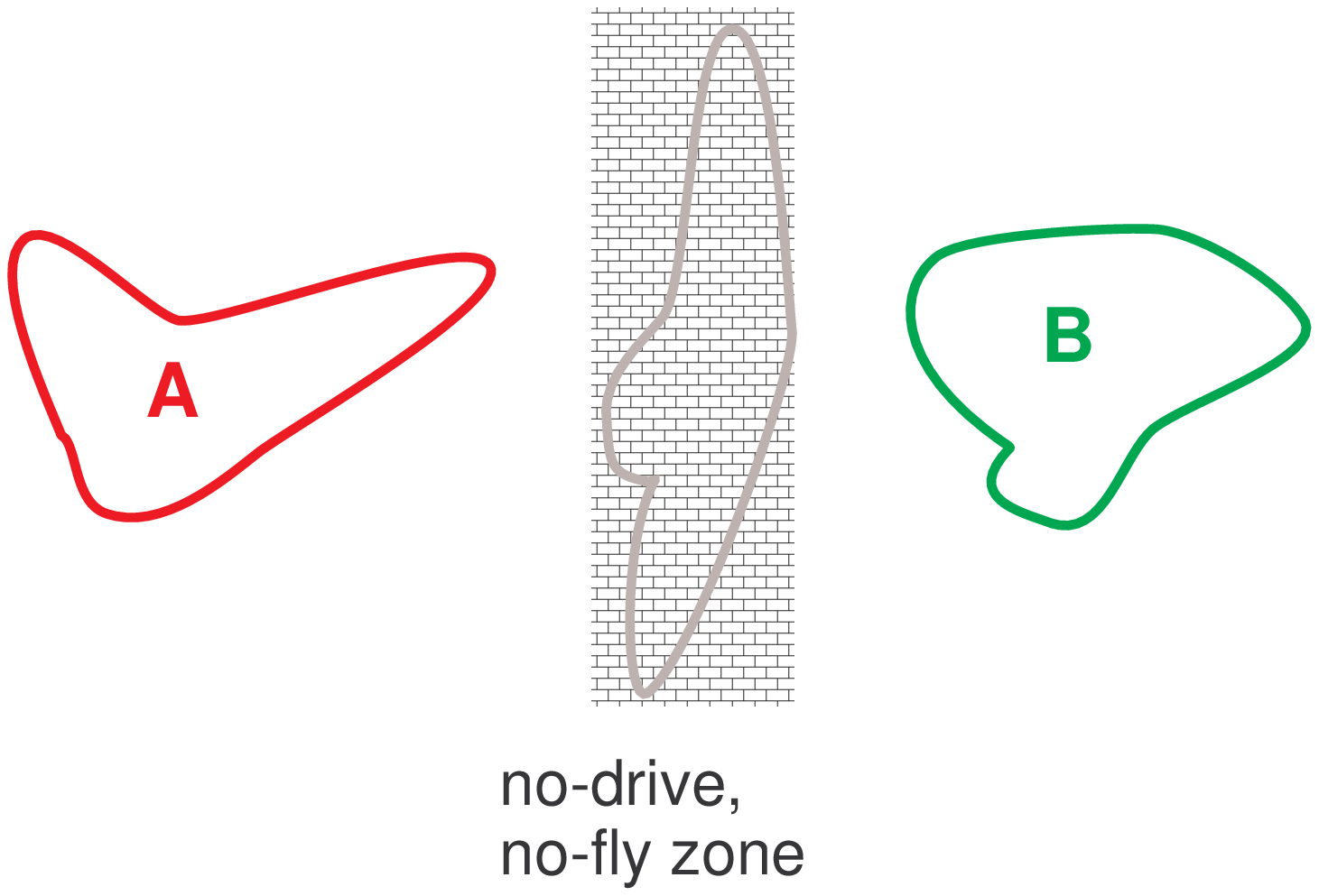}}
      }}
      \caption{\textsf{Illustrating some problems in defining distance between cities in the presence of obstacles.}}
    \label{fig:no-fly-zone}
   \end{figure}

\begin{figure}[h!]
      \centering
      \framebox{\parbox{0.8\textwidth}{\centering
      {\includegraphics[trim={1mm, 1mm, 1mm, 1mm}, clip, angle=-0, width = 0.75\textwidth]{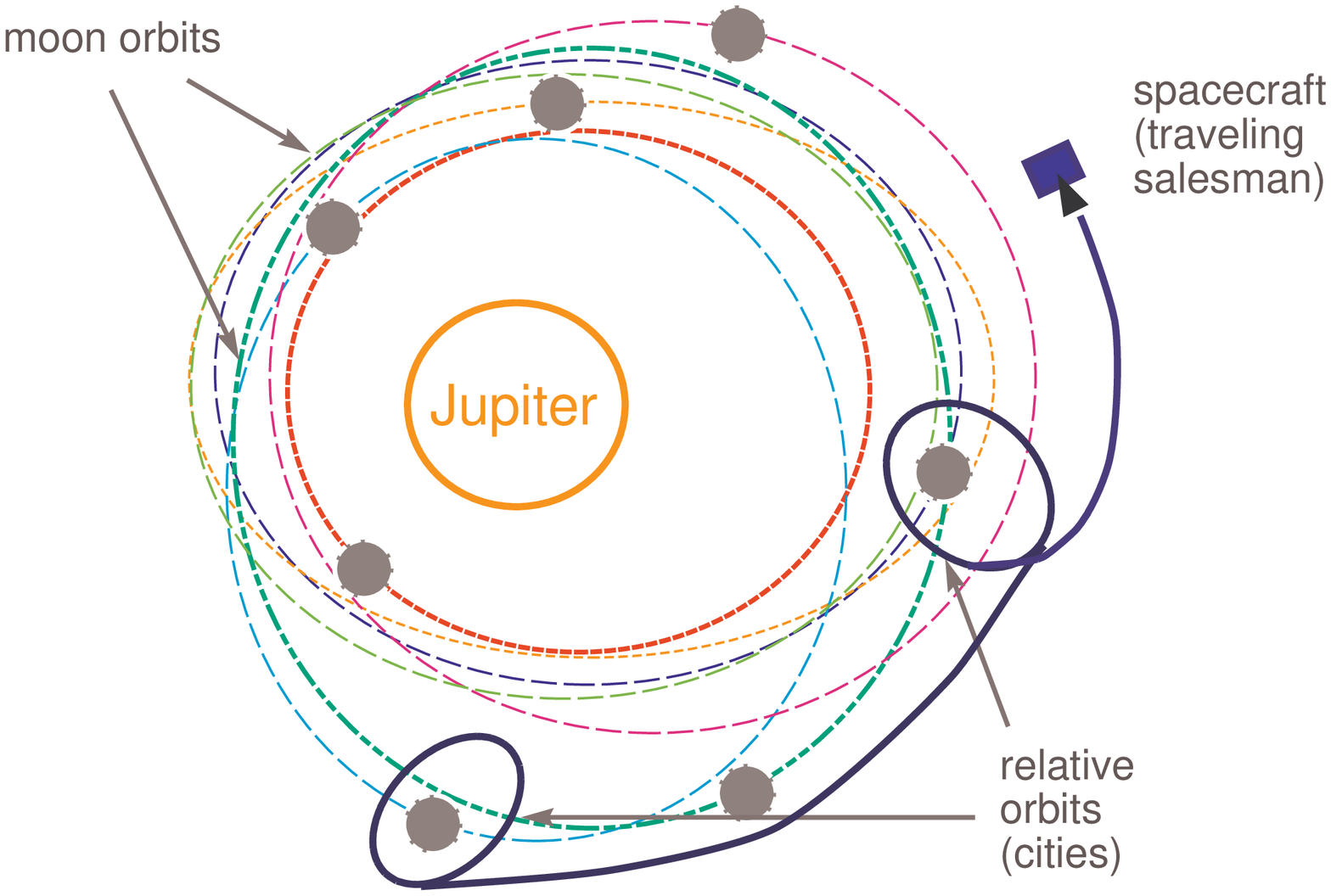}}
      }}
      \caption{\textsf{A rendition of a Jovian grand tour mission with only 7 of the 79 moons illustrated for clarity.  Figure is not to scale.}}
    \label{fig:Grand-tour}
   \end{figure}

Recently, we showed that a continuous-variable, nonlinear static optimization problem can be framed as a dynamic optimization problem \citep{ross-jcam-1}. In this theory, a generic point-to-set algorithmic map is defined in terms of a controllable continuous-time\footnote{Time is just a proxy for an independent variable.} trajectory, where, the decision variable is a continuous-time search vector.  Starting with this simple idea, many well-known algorithms, such as the gradient method and Newton's method, can be derived as optimal controllers over certain metric spaces \citep{ross-jcam-1}. If the control is set to the acceleration of a double-integrator model, then a similar theory\cite{ross-accel} generates accelerated optimization techniques such as Polyak's heavy ball method\cite{polyak64} and Nesterov's accelerated gradient method\cite{nesterov83}. In this paper, we further the theory of framing optimization problems in terms of an optimal control problem. More specifically, we use a menu of modifications to the traveling salesman problem (TSP) and its variants \citep{cook-2012,TSP-variants} to address a class of combinatorial optimization problems.  To motivate the new mathematical paradigm, consider the following questions:
%
%-----------------
\begin{question}
How do we define the distance between two sets if the vertices in a TSP graph are cities equipped with the power of the continuum?  Assume the cities to be disjoint sets; see Fig.~\ref{fig:2-City-Blobs}.
\end{question}
%-------------------
%%
%======================================================================================
%==========================================================================================
%%
A version of this question was posed in \cite{TSPN} more than two decades ago.  It continues to be an active research topic; see for example, \cite{TSPN-robotics} and \cite{wang2019}. To appreciate this question, consider a distance function defined by,
\begin{equation}\label{eq:d(AB)}
d(A,B) := \displaystyle\mathop\text{min}_{\bx \in A, \by \in B}d(\bx, \by) = \displaystyle\mathop\text{min}_{\bx \in A, \by \in B} \norm{\bx-\by}_2
\end{equation}
Besides the fact that $d(A,B)$ is not a metric, if \eqref{eq:d(AB)} is used to construct the edge weights in the TSP graph, the resulting solution may comprise disconnected segments because the entry and exit points for a city may not necessarily be connected. Adding a ``continuity segment''  \textit{a posteriori} does not generate an optimal solution as indicated by the three-city tour illustrated in Fig.~\ref{fig:3-City-Blobs}.
%
%======================================================================================
%==========================================================================================
%
This solution is clearly not optimal because the triangle tour shown in Fig.~\ref{fig:3-City-Blobs} is shorter.  In other words, an optimal tour is obtained by not using the shortest distance between two cities.
%---------------
\begin{remark}\label{rem:que-not-soln}
It is critically important to note that Question 1 is not centered on solving the problem of the type illustrated in Fig.~\ref{fig:3-City-Blobs}; rather, this question and others to follow, are focused more fundamentally on simply framing the problem mathematically.
\end{remark}
%---------------
\begin{remark}\label{rem:newTSP-objFun}
It is apparent by a cursory examination of Fig.~\ref{fig:3-City-Blobs} that the values of the arc weights are not independent of the path. That is, the objective function in the TSP must somehow account for the functional dependence of the sequence of cities in the computation of the distance between any two cites.
\end{remark}
%--------------
%--------------
\begin{question}
How do we define distances between two cities if the entry and exit points are constrained by some angle requirement?
\end{question}
%--------------
This question was first addressed in \cite{TSP-angle} for the case when the cities are points; the new question posed here pertains to the the additional issues when the cities are not points as in Fig.~\ref{fig:3-City-Blobs}.
%---------------
\begin{question}
How do we define distances between two cities in the presence of a no-drive, no-fly zone?  See Fig.~\ref{fig:no-fly-zone}.
\end{question}
%--------------
%
%======================================================================================
%==========================================================================================
%
This question was addressed in \cite{TSPFN} for the case of point-cities.  It is apparent that any difficulty encountered in answering Questions  2 and  3 is further amplified by the issues resulting from the discussions related to Question~1; see also Remarks \ref{rem:que-not-soln} and \ref{rem:newTSP-objFun}.
%------------------------
\begin{question}\label{que:moving-cities}
How do we define distances between neighborhoods that are in deterministic motion?
\end{question}
%-----------------------
This question is related to the one posed in \cite{DTSP} with respect to the dynamic TSP.  It remains a problem of ongoing interest; see, for example \cite{DTSP-2019}.
\begin{question}
Assuming we can answer the preceding questions, is distance a correct measure for a minimum-time TSP (with neighborhoods)? If not, what is the proper mathematical problem formulation for a minimum-time TSP?
\end{question}
%=========================
While limited versions of the aforementioned questions have been addressed in the literature as noted in the preceding paragraphs, the totality of Questions~1--5 appear in many practical and emerging mathematical problems that lie at the intersection of physics, operations research and engineering science. For example, the problem of touring the 79 moons of Jupiter \citep{witze} by a remote sensing spacecraft has all the elements of Questions 1--5. Relative orbits around the moons are the ``cities'' and the spacecraft is the ``traveling salesman;''  see Fig.~\ref{fig:Grand-tour}.
%
%======================================================================================
%==========================================================================================
%
The moons are in various non-circular orbits. The measure of ``distance'' (i.e., weights) is the amount of propellant it takes to transfer the spacecraft between two (moving) relative orbits.  Not all visits to the moons are valued equally by the science team.  The objective of a grand tour mission is to maximize the science return by orbiting around as many moons as possible under various constraints arising from the physics of gravitational motion, electromagnetic instrumentation, thermodynamics, electrical power, dollar cost and lifetime of the spacecraft. It is clear that modeling this optimization problem using the available constructs of a TSP \citep{cook-2012} is neither apparent nor easy.  In fact, the computation of the weights associated with the arcs of the graph (that represent the transfer trajectories) involve solving a constrained, nonlinear optimal control problem with variable endpoints \cite{ross-book,vinter}.  Consequently, even generating the data\footnote{To produce the TSP data, it is necessary to solve $\mathcal{O}(N^2)$ optimal control problems, where $N$ is the number of moons.} to define this problem as a standard TSP is a nontrivial task.

The main contribution of this paper is a new mathematical problem formulation for addressing a class of information-rich, operations-research-type mathematical problems such as the modified TSPs discussed in the preceding paragraphs.

\section{A New Mathematical Paradigm}
Throughout this paper we use the word functional in the sense of mathematical analysis: a mapping from a space of measurable functions to the field of real numbers.
%-----------------
\begin{definition}[$\mathcal{F}$-graph]
An $\mathcal{F}$-graph is a finite collection of functionals that constitute the arcs (edges) and vertices of a graph.
\end{definition}
%-----------------
Let $V^i: \text{dom}(V^i) \to \Real, i \in \mathbb{N}_+ $ and $E^k: \text{dom}(E^k) \to \Real, k \in \mathbb{N}_+$ be a finite collection of functionals, where $\text{dom}(\cdot)$ is the domain of $(\cdot)$.
Let $\mathcal{F}$ be an $\mathcal{F}$-graph whose vertices and arcs/edges are given by $V^i$ and $E^k$ respectively. From standard graph theory, a walk in $\mathcal{F}$ may be defined in terms of an alternating sequence of $V$- and $E$-functionals.  In order to perform evaluations in $\mathcal{F}$, it is necessary to define some new constructs.
%-----------------
\begin{definition}[$\mathcal{F}$-control]\label{def:F-control}
An $\mathcal{F}$-control is a sequence of functions
$$\langle \psi_0, \psi_1,  \ldots, \psi_n \rangle \quad n \in \mathbb{N}$$
where each $\psi_j, \ j = 0, \ldots, n$ is selected from the domain of $V^i$ or the domain of $E^k$.
\end{definition}
%----------------
\begin{remark}\label{rem:F-control}
An $\mathcal{F}$-control involves two simultaneous actions: selecting functions from the domain of the functionals that comprise $\mathcal{F}$, and ordering the selections in some sequence.  Thus, different selections ordered the same way is a different $\mathcal{F}$-control. Conversely, the same selection ordered differently is a also a different $\mathcal{F}$-control (provided the reordering is consistent with the selection process).
\end{remark}
%---------------
\begin{definition}[Control Walk]
A control walk (in $\mathcal{F}$) is an $\mathcal{F}$-control with the following property:
 $$ \psi_j \in \text{dom}(E^k)  \Rightarrow \psi_{j-1} \in \text{dom}(V^l) \ \text{and } \psi_{j+1} \in \text{dom}(V^m) $$
and $E^k$ is the arc/edge that joins $V^l$ to $V^m$.
\end{definition}
%----------------------------------
The preceding definition of a control walk requires a sequence of at least three functions. In order to complete this definition and accommodate certain special situations we define a trivial control walk as follows:
%-----------------
\begin{definition}[Trivial Control Walk]
A control walk $\omega_c$ is said to be trivial if:
\begin{enumerate}
\item $\omega_c = \langle \psi_0 \rangle$ and $\psi_0 \in \text{dom}(V^i)$ for some $i \in \mathbb{N}$; or,
\item $\omega_c = \langle \psi_0, \psi_1 \rangle$ and
$$\big(\psi_0 \in \text{dom}(V^i), \ \psi_1 \in \text{dom}(E^k) \big) \vee \big(\psi_0 \in \text{dom}(E^k), \ \psi_1 \in \text{dom}(V^i)\big)$$
for some $i, k$ in $\mathbb{N}_+$ and $E^k$ joins $V^i$ with itself or with some other vertex;
or
\item $\omega_c = \langle \psi_0, \psi_1, \psi_2 \rangle$ and
$$\psi_0 \in \text{dom}(E^k), \psi_1 \in \text{dom}(V^i), \psi_2 \in \text{dom}(E^m)$$
for some $i, k, m$ in $\mathbb{N}_+$, where $E^k$ and $E^m$ join $V^i$ with itself or with some other vertex.
\end{enumerate}
\end{definition}
%-----------------------------------------
\begin{remark}\label{rem:isomorp}
Let $G$ be a graph that is isomorphic to an $\mathcal{F}$-graph. Furthermore, let $G$ be such that its vertices and arcs are a selection of functions from the domains of the functionals the constitute the $\mathcal{F}$-graph. Then, by construction, there is an uncountably infinite set of isomorphic graphs $G$. A control walk may be interpreted as a walk in one of these isomorphic graphs.
\end{remark}
%------------------------------------
It follows from Remarks \ref{rem:F-control} and \ref{rem:isomorp} that a control walk involves the simultaneous action of choosing an isomorphism and a walk in the chosen isomorphic graph, $G$.

\begin{definition}[Objective Functional]
Let $\omega_c$ be a control walk defined over an $\mathcal{F}$-graph. An objective functional is a functional,
$\omega_c \mapsto \Real$,
defined over all possible control walks.
\end{definition}

\begin{definition}[Optimal Control Walk]
An optimal control walk is a control walk that optimizes an objective functional.
\end{definition}
%----------------------

%%----------------------
\begin{remark}
From Remark \ref{rem:F-control}, it follows that an optimal control walk jointly optimizes the selection of functions from the domains of the functionals that constitute an $\mathcal{F}$-graph as well as the walk itself. This feature of the optimal control walk directly addresses the comments in Remark \ref{rem:newTSP-objFun} in the context of Fig.~\ref{fig:3-City-Blobs}.
\end{remark}
%--------------------------------------

For the remainder of this paper, we will limit the discussions to a special type of an $\mathcal{F}$-graph defined over some measure space $\Lspace$.

%-----------------------------------------------
\begin{definition}[Label Space]
A measure space $\Lspace$ is called a label space if $\Lspace^i, i = 1, \ldots, |N_v| \in \mathbb{N}_+$ are a given finite collection of disjoint measurable subsets of $\Lspace$.
\end{definition}
%---------------------------------------

Suppose we are given a label space.  Let,
\begin{subequations}\label{eq:dom4Tg}
\begin{align}
\Lspace^a &:=\Lspace\setminus\displaystyle\bigcup_{i=1}^{|N_v|}\Lspace^i \\
\text{dom}(V^i) &:= \set{\Real \to \Lspace^i:\ \Real \to \Lspace^i \text{ is measurable}} \label{eq:domK4Tg}\\
\text{dom}(E^a) &:= \set{\Real \to \Lspace^a:\ \Real \to \Lspace^a \text{ is measurable}}\label{eq:domJ4Tg}
\end{align}
\end{subequations}
%
%---------------
\begin{definition}[$\mathcal{T}$-graph]
Let $E^{l,m}$ denote the arcs of an $\mathcal{F}$-graph with the property that $E^{l,m}$ joins $V^l$ to $V^m$ for all $l$ and $m$ in $N_v$.  Set $E^{l,m} = E^a$, where the domain of $E^a$ is given by \eqref{eq:domJ4Tg}. Set the domains of $V^i$ according to \eqref{eq:domK4Tg}. The resulting $\mathcal{F}$-graph is called a $\mathcal{T}$-graph.
\end{definition}
%-----------------------
\begin{definition}[Label Space Trajectory]
A label space trajectory is a measurable function $\bl(\cdot): \Real \supseteq [t_0, t_f] \ni t \mapsto \Lspace$ where $t_f - t_0 > 0$ is some non-zero time interval in $\Real$.
\end{definition}
%=====================================
\begin{proposition}\label{prop:1}
A label space trajectory is an $\mathcal{F}$-control for the $\mathcal{T}$-graph.
\end{proposition}
%------------------------
\proof{Proof:}
Let $\bl(\cdot): [t_0, t_f] \to \Lspace$ be a label space trajectory. Because $\Lspace^a$ and $\Lspace^i$ are all disjoint sets, we have,
\begin{equation}\label{eq:l_is_iore}
\bl(t) \in \Lspace^1 \vee \Lspace^2 \vee \ldots \vee\Lspace^{|N_v|}\vee \Lspace^a \text{ for a.a. } t \in [t_0, t_f]
\end{equation}
As a result, we can perform the following construction: Let,
$$t_0 < t_1, < \ldots, < t_{n} < t_{n+1}=t_f$$
be an increasing sequence of clock times such that for each $ j = 0, 1, \ldots, n,\ n \in \mathbb{N}$, $(t_{j+1} - t_j)$ is the maximum time duration for which we have,
\begin{equation}\label{eq:lseg_is_iore}
\bl(t) \in \left\{
                \begin{array}{ll}
                  \Lspace^i \text{ for some } i \in N_v & \hbox{and for a.a. }  t \in (t_j, t_{j+1})  \\
& \text{or}\\
                  \Lspace^a & \hbox{for a.a. } t \in (t_j, t_{j+1})
                \end{array}
              \right.
\end{equation}
The feasibility of constructing \eqref{eq:lseg_is_iore} follows from \eqref{eq:l_is_iore}.
Let,
%%
%$$\psi_j: (t_{j} - t_{j-1}) \ni t \mapsto \bl(t), \ j = 1, \ldots, n$$
%%
%
\begin{equation}\label{eq:psi-def-restrict}
\psi_j:= \bl(\cdot)\big|_{(t_{j} , t_{j+1})}, \ j = 0, 1, \ldots, n
\end{equation}
be the restrictions of $\bl(\cdot)$ to $(t_{j} , t_{j+1})$ for $j = 0, 1, \ldots, n$.  Then, by \eqref{eq:lseg_is_iore} and \eqref{eq:psi-def-restrict}, we have,
%
%%
%\begin{equation}\label{eq:Impsi}
%\psi_j(t) \in \left\{
%                \begin{array}{ll}
%                  \Lspace^i & \hbox{for a.a. }  t \in (t_j, t_{j+1})  \\
%& \text{or}\\
%                  \Lspace^a & \hbox{for a.a. } t \in (t_j, t_{j+1})
%                \end{array}
%              \right. \qquad\forall\ j = 0, 1, \ldots, n
%\end{equation}
%%
%
\begin{equation}\label{eq:Impsi}
\text{Im}(\psi_j) \subseteq \Lspace^i \vee \Lspace^a \text{ for some } i \in N_v \text{ and } \forall\ j = 0, 1, \ldots, n
\end{equation}
Hence from \eqref{eq:dom4Tg}, it follows that,
$$\psi_j \in \text{dom}(V^i) \vee \text{dom}(E^a), \text{ for some } i \in N_v$$
Consequently the sequence $\langle \psi_0, \ldots, \psi_n \rangle$ is an $\mathcal{F}$-control for the $\mathcal{T}$-graph.
%
%Let $[t_0, t_f] \ni t \mapsto \psi(t)$ be defined by $\psi_j(t)$ if $t \in (t_j, t_{j+1})$.  Then, from \eqref{eq:lseg_is_iore} and \eqref{eq:psi-def-restrict}, it follows that $\bl(t) = \psi(t)$ for $a.a.\ t \in [t_0, t_f]$.
%
%$,\ \bl(t) = \psi_j(t) \text{ for some } j = 0, 1, \ldots, n  $
%
% $\bl(\cdot)$ is an $\mathcal{F}$-control for the $\mathcal{T}$-graph.
%
\endproof
%======================================
%--------------------
\begin{definition}[Trivial Label Space Trajectory]
A label space trajectory is said to be trivial if $\textrm{Im}(\bl(\cdot)) \subseteq \Lspace^a$.
\end{definition}
%==============================
\begin{theorem}
A nontrivial, continuous label space trajectory generates a control walk in a $\mathcal{T}$-graph.
\end{theorem}
%---------------------------------------
\proof{Proof:}
Let $\psi_j, j = 0, \ldots, n \in \mathbb{N}$ be the restrictions of $\bl(\cdot)$ as defined by \eqref{eq:psi-def-restrict}. By Proposition \ref{prop:1}, $\langle \psi_0, \psi_1, \ldots, \psi_n \rangle$ is an $\mathcal{F}$-control; hence, any given $\psi_j$ is either in the domain of $E^a$, or of $V^i$, for some $i \in N_v$.  The rest of the proof is broken down to three cases:

\noindent\textbf{Case(a)}: $n \ge 2$ and $\psi_j \in \text{dom}(E^a)$ for some $0 < j < n$. \\
From \eqref{eq:lseg_is_iore} and \eqref{eq:psi-def-restrict} we get the following conditions:
\begin{enumerate}
\item $\psi_{j-1} \in \text{dom}(V^l)$ for some $l \in N_v$
\item $\psi_{j+1} \in \text{dom}(V^m)$ for some $m \in N_v$
\end{enumerate}
Hence, $\langle \psi_{j-1}, \psi_j, \psi_{j+1} \rangle$ is a control subwalk.

\noindent\textbf{Case(b)}: $n \ge 2$ and $\psi_j \in \text{dom}(V^i)$ for some $0 < j < n $ and some $i \in N_v$.\\
By continuity of $\bl(\cdot)$, we have,
\begin{equation}\label{eq:l-psi-continuity}
\lim_{t \uparrow t_{k+1} }\psi_k(t) = \lim_{t \downarrow t_{k+1} } \psi_{k+1}(t), \quad k = 0, \ldots, n
\end{equation}
Setting $k = j-1$ and $k = j$ in \eqref{eq:l-psi-continuity} we get the two continuity conditions,
\begin{subequations}
\begin{align}
\lim_{t \uparrow t_{j} }\psi_{j-1}(t) &= \lim_{t \downarrow t_{j} } \psi_{j}(t)\\
\lim_{t \uparrow t_{j+1} }\psi_j(t) &= \lim_{t \downarrow t_{j+1} } \psi_{j+1}(t)
\end{align}
\end{subequations}
Hence, from \eqref{eq:domJ4Tg}, \eqref{eq:domK4Tg}, \eqref{eq:lseg_is_iore} and \eqref{eq:psi-def-restrict}, we have,
\begin{equation}
\psi_{j-1} \in \text{dom}(E^a)  \text{ and }  \psi_{j+1} \in \text{dom}(E^a)
\end{equation}
If $n=2$ ($\Rightarrow j=1$), we get a trivial control walk.  If $ j >1$, by using the same arguments as in Case(a), we get $\psi_{j-2} \in \text{dom}(V^k) \text{ for some } l \in N_v \text{ and }  \psi_{j+2} \in \text{dom}(V^m) \text{ for some } m \in N_v$; hence, $\langle \psi_{j-2}, \psi_{j-1}, \psi_j, \psi_{j+1}, \psi_{j+2} \rangle$ is a control subwalk.

\noindent\textbf{Case(c)}: $n<2$.\\
By similar arguments as in Case (a) and (b), it is straightforward to show that $\langle \psi_0, \psi_1 \rangle$ is trivial control walk.

From Cases (a), (b) and (c), it follows that $\langle \psi_0, \psi_1, \ldots, \psi_n \rangle, n \in \mathbb{N}$ is a control walk.
\endproof
%=======================================

%\section{A $\mathcal{T}$-Graph Formulation of a TSP in Label Space}
\section{Two $\mathcal{T}$-Graph-Based Formulations of a TSP}

Let $\Lspace$ be a finite dimensional normed space.  If we set $\Lspace^i$ as the ``cities'' in $\Lspace$, then a TSP and its many variants may be framed in terms of finding nontrivial optimal label-space trajectories.  To illustrate the theoretical simplicity of our approach, we allow $\Lspace^i$ to be time dependent (i.e., in deterministic motion) so that some answers to the questions posed in Section 1 are readily apparent.
%Furthermore, we use a Kronecker indicator function to generate the key construct of a restriction operation used in Proposition 1 (Cf. \eqref{eq:lseg_is_iore}):
%
%
\begin{definition}[Atomic Return Function]
For a fixed time $t$, an atomic return function $R_a : (\bl, \Lspace^i(t)) \mapsto \Real$ is defined by,
%for a fixed time $t$
% by the following two properties:
%%
%\begin{subequations}\label{eq:return-def}
%\begin{align}
%&& R_a(\bl, \Lspace^i(t)) &= 0         &\text{if } \bl \not\in \Lspace^i(t)  \\
%&& \int_{\Lspace^i(t)} R_a(\bl, \Lspace^i(t))\, d\mu  &\neq 0  &\text{if } \bl \in \Lspace^i(t)
%\end{align}
%\end{subequations}
%%
%
\begin{equation}\label{eq:atomic-ret-fun}
R_a(\bl, \Lspace^i(t))  \left\{
                         \begin{array}{lll}
                         \neq 0 & \hbox{if} \ \bl \in \Lspace^i(t) \\
                           = 0 & \hbox{otherwise}
                         \end{array}
                       \right.
\end{equation}
%
%where,
%%
%\begin{equation}\label{eq:return-supp}
%\int_{\Lspace^i(t)} R_a(\bl, \Lspace^i(t))\, d\mu  \ne 0
%\end{equation}
%%
%$\mu$ is a measure over $\Lspace^i(t)$.
\end{definition}
%-------------------
%-------------------
\begin{definition}[Atomic Return Functional]\label{def:atomic-ret-func}
An atomic return functional $R^i$ is defined by,
\begin{equation}\label{eq:return-functional-def}
R^i[\bl(\cdot), t_0, t_f] := \int_{t_0}^{t_f} R_a\big(\bl(t), \Lspace^i(t)\big) dt
\end{equation}
where $t_f - t_0 > 0$ is the time horizon of interest.
\end{definition}
%----------------------
Atomic return functionals may be used as vertex functionals whenever $\bl(t) \in \Lspace^i(t)$.  Whether or not a vertex functional is defined for a given problem, we define the Kronecker indicator function,
\begin{equation}\label{eq:kron-indicator}
\mathcal{I} (\bl, \Lspace^i(t)) := \left\{
                         \begin{array}{ll}
                           1 & \hbox{if } \bl \in \Lspace^i(t) \\
                           0 & \hbox{otherwise}
                         \end{array}
                       \right.
\end{equation}
as a fundamental return function.  Equation \eqref{eq:kron-indicator} thus generates by way of \eqref{eq:return-functional-def} a fundamental return functional:
%----------------------------
\begin{definition}[Time-on-Task Functional] A time-on-task functional is defined by,
\begin{equation}\label{eq:Ti-def}
 T^i[\bl(\cdot), t_0, t_f]:= \displaystyle \int_{t_0}^{t_f} \mathcal{I} (\bl(t), \Lspace^i(t))\, dt
\end{equation}
\end{definition}
%---------------------------
Because \eqref{eq:Ti-def} generates a dwell time over vertex $i$, we define a visit in the context of a walk in the following manner:
%
%-------------------
\begin{definition}[Vertex Visit]
The vertex $i$ is said to have been visited in $[t_0, t_f]$ if
\begin{equation}\label{eq:visit-def}
T^i[\bl(\cdot), t_0, t_f] \ne 0
\end{equation}
Correspondingly, we say the vertex $i$ has not been visited if $T^i[\bl(\cdot), t_0, t_f] = 0$.
\end{definition}
%-------------------
%
In the context of a $\mathcal{T}$-graph formulation of a basic TSP,\footnote{By a basic TSP, we mean a problem that does not come with any additional qualifiers.} a vertex functional $V^i$ assigns a value of one for a visit, zero otherwise.  Furthermore, an arc functional $E^a$ is any functional that assigns a numerical value (of ``distance'') for segments of $t \mapsto \bl$ that are not associated with a vertex. If a visit is preceded and followed by an arc functional with no other visits of a vertex in between (including itself), we say vertex $i$ has been visited once.  In this context, we say $\bl(\cdot)$ is Hamiltonian if all vertices are visited exactly once.

\subsection{A Derivative-Based Formulation of a TSP}

%=================================
\begin{lemma}
Let $D^i[\bl(\cdot), t_0, t_f]$ be the functional defined by,
\begin{equation}
D^i[\bl(\cdot), t_0, t_f] := \int_{t_0}^{t_f}\abs{d_t\mathcal{I} \left(\bl(t), \Lspace^i(t)\right)}dt
\end{equation}
where, $d_t$ denotes the distributional derivative with respect to $t$. Let $\bl(\cdot): [t_0, t_f] \to \Lspace$ be a continuous label space trajectory such that $T^i[\bl(\cdot), t_0, t_f] \ne 0$. Furthermore, let $\bl(t_0) \not\in \Lspace^i(t_0)$ and $\bl(t_f) \not\in \Lspace^i(t_f)$. Then,
\begin{equation}\label{eq:Di=even}
D^i[\bl(\cdot), t_0, t_f]  \in 2\mathbb{N}_+
\end{equation}
%
%%
%\begin{equation}\label{eq:no-lingering}
% \int_{\partial\Lspace^i(t)}\bl(t)\, dt =0
%\end{equation}
%%
%where, $\partial\Lspace^i(t)$ is the boundary of $\Lspace^i(t)$. Then,
%%
%\begin{equation}\label{eq:Di-in-N}
%D^i[\bl(\cdot), t_0, t_f]  \in \mathbb{N}
%\end{equation}
%%
%
\end{lemma}
%------------------------------------
\proof{Proof:}
Consider the time intervals $\Delta t_{out}$ and $\Delta t_{in}$ defined by,
\begin{equation}
\Delta t_{out}:= \set{t \in [t_0, t_f]:\ \mathcal{I}(\bl(t), \Lspace^i(t) = 0 } \quad \Delta t_{in}:= \set{t \in [t_0, t_f]:\ \mathcal{I}(\bl(t), \Lspace^i(t) = 1 }
\end{equation}
By  assumption $T^i[\bl(\cdot), t_0, t_f] \ne 0$; hence, we have $\mu\left(\Delta t_{in}\right) > 0 $.  Likewise, we have $\mu\left(\Delta t_{out}\right) > 0 $.  Hence, we can write,
$$ \frac{d\mathcal{I} \left(\bl(t), \Lspace^i(t)\right)}{dt} = 0 \quad\forall\ t \in \set{\text{int}\left(\Delta t_{out}\right)} \cup \set{\text{int}\left(\Delta t_{in}\right)}$$
where, int$(\cdot)$ denotes the interior of the set $(\cdot)$. Let
$ t_x \in \partial\Delta t_{int} $
where $\partial(\cdot)$ denotes the boundary of the set $(\cdot)$.  Then, the distributional derivative of the function $t \mapsto   \mathcal{I} \left(\bl(t), \Lspace^i(t)\right)$ may be written as,
\begin{equation}\label{eq:diff-I}
 d_t\mathcal{I} \left(\bl(t), \Lspace^i(t)\right)=                                                                        -\delta(t-t_x) \vee \delta(t-t_x)
\end{equation}
%
%%
%\begin{equation}\label{eq:diff-I-cases}
% d_t\mathcal{I} \left(\bl(t), \Lspace^i(t)\right)\big\vert_{t_x} =\left\{
%                                                                      \begin{array}{ll}
%                                                                        -\delta(t-t_x) \vee \delta(t-t_x) & \hbox{if } \bl(t_x) \in \partial\Lspace^i(t_x) \\
%                                                                        0 & \hbox{if } \bl(t_x) \not\in \partial\Lspace^i(t_x)
%                                                                      \end{array}
%                                                                    \right.
%\end{equation}
%%
where, $\delta(t-t_x)$ is the Dirac delta function centered at $t=t_x$. Hence,
\begin{equation}\label{eq:diff-I-abs=delta}
\abs{d_t\mathcal{I} \left(\bl(t), \Lspace^i(t)\right)} = \delta(t-t_x)
\end{equation}
Because $\bl(t_0)$ and $\bl(t_f)$ are not in $\Lspace^i(t_0)$ and $\Lspace^i(t_f)$ respectively, integrating \eqref{eq:diff-I-abs=delta} we get,
$$ \int_{t_0}^{t_f}\abs{d_t\mathcal{I} \left(\bl(t), \Lspace^i(t)\right)}dt  \in 2\mathbb{N}_+ $$
where, we have used the fact that the integral of a delta function is unity.
\endproof
%======================================================
\begin{corollary}
If $\bl(t_0) \in \Lspace^i(t_0)$ and $\bl(t_f) \in \Lspace^i(t_f)$, then \eqref{eq:Di=even} generalizes to $D^i[\bl(\cdot), t_0, t_f]  \in 2\mathbb{N}$. If $\bl(t_0) \in \Lspace^i(t_0)$ or $\bl(t_f) \in \Lspace^i(t_f)$, then the statement of Lemma 1 may be further generalized to $D^i[\bl(\cdot), t_0, t_f]  \in \mathbb{N}$.
\end{corollary}
%====================================

%=================================
\begin{theorem}
Let $\bl(\cdot): [t_0, t_f] \to \Lspace$ be a continuous label space trajectory. Then, $\bl(\cdot)$ is Hamiltonian if and only if
\begin{equation}\label{eq:Di=2}
 D^i[\bl(\cdot), t_0, t_f] = 2 \quad \forall\ i \in N_v
\end{equation}
\end{theorem}
%------------------------------------
\proof{Proof:}
If $\bl(\cdot)$ is Hamiltonian, then every vertex has been visited just once; hence, \eqref{eq:Di=2} follows from Lemma 1.  If \eqref{eq:Di=2} holds, then every vertex has been visited just once in $[t_0, t_f]$ (with the understanding of a visit given by \eqref{eq:visit-def}).
\endproof
%======================================================
Although continuity of a label space trajectory is sufficiently smooth to generate a Hamiltonian cycle (Cf. Theorem 1), we now consider the space of absolutely continuous functions
%$AC([t_0, t_f), \Lspace)$,
for the convenience of using the derivative of $t \mapsto \bl$ to define arc lengths. To this end, let $\dot\bl(t)$ denote the time derivative of $\bl(t)$; then,  a distance functional may be defined according to,
\begin{equation}\label{eq:J=dist}
J_{dist}[\bl(\cdot), t_0, t_f] := \displaystyle \int_{t_0}^{t_f} \norm{\dot\bl(t)} dt
\end{equation}
The integrand in \eqref{eq:J=dist} is any finite-dimensional norm. If the two-norm is used, then the numerical value of $J_{dist}$ is consistent with the notion of Euclidean distance illustrated in Fig.~\ref{fig:3-City-Blobs}.  Combining \eqref{eq:J=dist} with Theorem 2, a shortest distance TSP can be formulated as,
%==========================
\begin{eqnarray}
&\big(\text{$D$-TSP}\big) \left\{
\begin{array}{lll}
\text{Minimize}  & J_{dist}[\bl(\cdot), t_0, t_f] := \displaystyle \int_{t_0}^{t_f} \norm{\dot\bl(t)} dt \\[1.1em]
\text{Subject to }
& \displaystyle D^i[\bl(\cdot), t_0, t_f] = 2 \quad \forall\ i \in N_v\\[0.7em]
& \bl(t_f) = \bl(t_0)
\end{array} \right.& \label{eq:prob-TSP-classic}
\end{eqnarray}
%========================================
%
The last constraint equation in \eqref{eq:prob-TSP-classic} simply ensures that the resulting label-space trajectory is a closed control-walk.

In comparing it with the various discrete variable optimization formulations of a TSP \citep{Dantzig-54,Flood,MTZ},  it is apparent that \eqref{eq:prob-TSP-classic} contains  no explicit subtour-type elimination constraints.  This is because the (absolute) continuity of the label space trajectory ensures (by Theorem 1) that the closed control-walk generated by \eqref{eq:prob-TSP-classic} is a single Hamilton cycle.

\subsection{An Integral-Based Formulation of a TSP}

A TSP may also be formulated (in the sense of a $\mathcal{T}$-graph) using the time-on-task functional  to construct a new indicator-type functional.
%
%----------------------------
\begin{definition}[Control-Walk Indicator Functional] A control-walk indicator functional is defined by,
\begin{equation}\label{eq:Wi-def}
W^i[\bl(\cdot), t_0, t_f] :=  \left\{
                         \begin{array}{lll}
                           1 & \hbox{if } \  T^i[\bl(\cdot), t_0, t_f] \ne 0 \\[1.1em]
                           0 & \hbox{if } \  T^i[\bl(\cdot), t_0, t_f] = 0
                         \end{array}
                       \right.
\end{equation}
\end{definition}
%---------------------------
%
Using $W^i$ as an indicator of a vertex visit, we arrive at an alternative formulation of a TSP:
%
%==========================
\begin{eqnarray}
&\big(\text{$I$-TSP}\big) \left\{
\begin{array}{lll}
\text{Minimize}  & J_{dist}[\bl(\cdot), t_0, t_f] := \displaystyle \int_{t_0}^{t_f} \norm{\dot\bl(t)} dt \\[1.1em]
\text{Subject to }
& \displaystyle W^i[\bl(\cdot), t_0, t_f] = 1 \quad \forall\ i \in N_v\\[0.7em]
& \bl(t_f) = \bl(t_0)
\end{array} \right.& \label{eq:prob-TSP-classic-I}
\end{eqnarray}
%=======================================
Equations \eqref{eq:prob-TSP-classic-I} and \eqref{eq:prob-TSP-classic} are identical except for the imposition of degree constraints.  Formulation ($I$-TSP) requires that the vertex be visited at least once.
%
%=======================
\begin{table}[h!]
\begin{center}
\begin{tabular}{|l|p{0.33\linewidth}|p{0.12\linewidth}|r|} %{|p{0.3cm}|p{3.5cm}|p{1.75cm}|r|}
  \hline \hline
No & Entity  & Discrete Optimization & $\mathcal{T}$-Graph Formulations \\
\hline
1 &  optimization variable & discrete & continuous function \\
2 & cost matrix; i.e., explicit arc/edge weights & required &  not required \\
3 & explicit degree constraints & required & required in \eqref{eq:prob-TSP-classic}; not in \eqref{eq:prob-TSP-classic-I} \\
4 & explicit subtour elimination constraints & required & not required\\
  \hline
\end{tabular}\\[0.5em]
\caption{A discriminating comparison between a discrete-variable and $\mathcal{T}$-graph formulations of a TSP }
\end{center}
\end{table}
%=======================
A comparison between the $\mathcal{T}$-graph and discrete-variable problem formulations is summarized in Table 1.  As noted in the second row of Table~1, the $\mathcal{T}$-graph formulations do not require the explicit construction of the arc/edge weights; i.e., the cost matrix for the computation of the objective function in the discrete-variable formulation.  This data is ``generated'' simultaneously via \eqref{eq:J=dist} as part of the process of solving the problem.  That is, in terms of the concepts introduced in Section 2, the $\mathcal{T}$-graph formulation incorporates the simultaneous selection of the function sequence and the functions themselves from the domains of the vertex and arc functionals.
%-------------
\begin{remark}
The $\mathcal{T}$-graph formulations presented in this section are not transformations of the well-established discrete-optimization models of TSPs; rather, they are a realization of a fundamentally new domain of analysis.
\end{remark}
%----------------

\section{Sample $\mathcal{T}$-Graph Formulations of Several Variants of a TSP}
Because of the large number of variants of a TSP, we limit the scope of this section to a small sample of selective problems to illustrate the new modeling framework.
%In formulating some of the variants of a TSP, it is convenient to use the class of continuous functions $C([t_0, t_f), \Lspace)$, in addition to the space of absolutely continuous functions used in Section 3.

\subsection{An Orienteering Problem}
Let $\sigma^i > 0$ be the score associated with each $\Lspace^i(t)$.  Define a score functional according to:
$$S^i[\bl(\cdot), t_0, t_f] := \sigma^i W^i[\bl(\cdot), t_0, t_f]  $$
where, $W^i$ is given by \eqref{eq:Wi-def}. An orienteering problem (OP) may now be defined as,
%==========================
\begin{eqnarray}
&\big(\text{OP}\big) \left\{
\begin{array}{lll}
\text{Maximize}  & \displaystyle \sum_{i \in N_v} S^i[\bl(\cdot), t_0, t_f] \\[1.1em]
\text{Subject to}
& \displaystyle \int_{t_0}^{t_f} 1\, dt \le  t_{max}  \\[0.7em]
& \bl(t_0) \in \Lspace^0 \\
& \bl(t_f) \in \Lspace^{|N_v|}
\end{array} \right.& \label{eq:prob-OP-classic}
\end{eqnarray}
%========================================
%
The payoff functional is the sum of the score functionals.  The maximum allowable time is $t_{max} > 0$.  The time constraint in \eqref{eq:prob-OP-classic} is written in terms of the integral of one to merely illustrate the fact that the resource constraint is a functional. The domain of $\bl(\cdot)$  in \eqref{eq:prob-OP-classic} is the space of continuous functions in accordance with Theorem 1.

\subsection{An Orienteering Problem With Neighborhoods}
Problem (OP) given by \eqref{eq:prob-OP-classic} was motivated by the usual orienteering problem discussed in the literature \citep{OP-survey}. Using the atomic return function concept introduced in \eqref{eq:atomic-ret-fun}, we may now score a traversal through a neighborhood based on the values of the atomic return functional (Cf.~\eqref{eq:return-functional-def}).  Although there are many ways to score a traversal, we illustrate a problem formulation based on the following construction,
\begin{equation}\label{eq:S=nbhd-step}
S^i_{nbd}[\bl(\cdot), t_0, t_f] :=   \left\{
                         \begin{array}{lll}
                           \sigma^i & \hbox{if} \ R^i[\bl(\cdot), t_0, t_f] \ge r^i\\
                           0 & \hbox{if} \ R^i[\bl(\cdot), t_0, t_f] < r^i
                         \end{array}
                       \right.
\end{equation}
In \eqref{eq:S=nbhd-step}, $r^i > 0$ is a required ``revenue'' from a visit to $\Lspace^i(t)$. If a label-space traversal $t \mapsto \bl$ across the city does not generate a revenue of at least $r^i$ as measured by the return functional $R^i$, then the salesman gets no commission (i.e., zero score). On the other hand, if the salesman performs a judicious travel through the city to generate a revenue of at least $r^i$, then, he is rewarded by $\sigma^i > 0$.  The salesman makes no extra commission for generating a revenue greater than $r^i$ at $\Lspace^i$; i.e., he is encouraged to expand the market by visiting a different city. This situation arises in astronomy \citep{cook-2012} where a telescope is required to scan a portion of the sky to collect a specific frequency of electromagnetic (EM) radiation \citep{TSP-telescope}.  Scientific value is generated based on the integration of EM collects. A critical amount of EM collects ($r^i$) is necessary to perform useful science. No extra science is generated once the targeted amount $r^i$ is reached; hence, the telescope is awarded $\sigma^i$ for performing a task successfully with no ``extra'' credit. Thus, the orienteering problem with neighborhoods may be defined according to,
%
%==========================
\begin{eqnarray}
&\big(\text{OP-nbd}\big) \left\{
\begin{array}{lll}
\text{Maximize}  & \displaystyle \sum_{i \in N_v} S^i_{nbd}[\bl(\cdot), t_0, t_f] \\[1.1em]
\text{Subject to }
& \displaystyle C^k[\bl(\cdot), t_0,t_f] \le  C^k_{max} \quad \forall\ k \in N_C  \\[0.7em]
& \bl(t_0) \in \Lspace^0 \\
& \bl(t_f) \in \Lspace^{|N_v|}
\end{array} \right.& \label{eq:prob-OP-nbhd}
\end{eqnarray}
%========================================
%
In \eqref{eq:prob-OP-nbhd}, the constraints $C^k[\bl(\cdot), t_0, t_f] \le  C^k_{max},\ \forall\ k \in N_C$ are $N_C \subset \mathbb{N}$ generic resource constraint, where $C^k_{max} > 0$ is the maximum ``capacity'' associated with the $k$-th resource.

\subsection{A Fast TSP}
Minimizing distance traveled does not necessarily equate to minimum time; this fact has been known since Bernoulli posed his famous Brachistochrone problem as a mathematical challenge in the year 1696 \cite{ross-book,vinter}. In framing a more interesting minimum-time TSP, we let $\tau^i > 0$ be an additional constraint of a minimum required dwell-time over $\Lspace^i(t)$. In this case, a fast (i.e., minimum-time) TSP may be framed as follows:
%==========================
\begin{eqnarray}
&\big(\text{fastTSP}\big) \left\{
\begin{array}{lll}
\text{Minimize}   &J_{time}[\bl(\cdot), t_0, t_f]:= \displaystyle \int_{t_0}^{t_f} 1\, dt \\[1.1em]
\text{Subject to}
& T^i[\bl(\cdot), t_0, t_f] \ge  \tau^i \qquad \forall\ i \in N_v\\[0.7em]
& \bl(t_f) = \bl(t_0)
\end{array} \right.& \label{eq:prob-TSP-minTime}
\end{eqnarray}
%========================================
%
For the same reason as the formulation of the constraint in \eqref{eq:prob-OP-classic}, the cost function in \eqref{eq:prob-TSP-minTime} is written as an integral to emphasize the fact that the travel-time is a functional.  The functional $T^i$ in \eqref{eq:prob-TSP-minTime} is given by \eqref{eq:Ti-def}.  Furthermore, while the dwell-time constraint in \eqref{eq:prob-TSP-minTime} may, in principle, be written as an equality, $ T^i[\bl(\cdot), t_0, t_f] =  \tau^i $, the advantage of an inequality is that it allows the solution to satisfy the stipulated requirement without incurring a penalty in performance for ``navigating through the city'' to find an optimal exit point at an optimal exit time.

\subsection{Dynamic TSP with Time Windows}
Because $\Lspace^i(t), i \in N_v$ are ``moving cities,'' dynamic TSPs are are explicitly incorporated in all of the preceding problem formulations. These formulations have also implicitly incorporated time windows because of the following argument: Let $t_a^i \in [t_0, t_f]$ and  $t_b^i \in [t_0, t_f]$ be two given clock-times associated with each $\Lspace^i(t)$  with $t_a^i < t_b^i$. Consider the augmented set defined by,
\begin{equation}\label{eq:Lspace-i-window}
\Lspace^i_{aug}(t):= \Lspace^i(t) \times [t_a^i, t_b^i]
\end{equation}
Then, it is clear that the space defined by,
$$\Lspace_{aug}(t):= \Lspace(t) \times [t_0, t_f]$$
is a label space, whose disjoint sets are defined by \eqref{eq:Lspace-i-window}.  Hence, by defining a new label-space variable $\bl_{aug} := (\bl, t) $, it is apparent that no additional theoretical developments are necessary to incorporate time-window constraints in the proposed $\mathcal{T}$-graph formulations.  Furthermore, note that the time window in \eqref{eq:Lspace-i-window} may itself also vary with respect to time.

In view of economic considerations that go beyond distance and time, a significantly greater degree of flexibility in modeling can be obtained by transforming the preceding $\mathcal{T}$-graph ``variational'' formulations to their optimal-control versions.

\section{A $\mathcal{T}$-Graph Optimal Control Framework}

By limiting the space of allowable label space trajectories to the space of absolutely continuous functions, it is possible to frame a significantly richer class of traveling-salesman-type problems. To develop this transformed framework, we first set the derivative of $\bl(t)$ as a candidate optimization variable,
\begin{equation}\label{eq:ldot=w}
\dot\bl(t) = \bw(t)
\end{equation}
where, $\bw \in \real{N_l}, N_l \in \mathbb{N}$ is the label-space tangent control variable.  Thus, for example, \eqref{eq:prob-TSP-minTime}  transforms according to,
%==========================
\begin{eqnarray}
&\big(\text{fastTSP: OC-ver}\big) \left\{
\begin{array}{lll}
\displaystyle\mathop\text{Minimize }   & \displaystyle J_{time}[\bl(\cdot), \bw(\cdot), t_0, t_f] := t_f - t_0 \\[1.1em]
\text{Subject to }
& \dot\bl(t) = \bw(t)\\[0.7em]
& T^i[\bl(\cdot), t_0, t_f] \ge  \tau^i \qquad \forall\ i \in N_v\\[0.7em]
& \bl(t_f) = \bl(t_0)
\end{array} \right.& \label{eq:prob-TSP-minTime-ocp-basic}
\end{eqnarray}
%========================================
%
%The vertex functional constraint is an isoperimetric-type constraint that is well-known in optimal control theory and goes back to famous problem of Queen Dido in the Badlands [CLARKE].
%Other problems discussed in Section 3 can be similarly transformed.
%
In comparing it with \eqref{eq:prob-TSP-minTime}, note that \eqref{eq:prob-TSP-minTime-ocp-basic} contains an additional decision variable $\bw(\cdot)$ and an additional constraint given by \eqref{eq:ldot=w}.  This aspect of the transformation may be used advantageously; for example,
if the objective functional in \eqref{eq:prob-TSP-minTime-ocp-basic} is replaced by,
$$J_{dist}[\bl(\cdot), \bw(\cdot), t_0, t_f] := \displaystyle \int_{t_0}^{t_f} \norm{\bw(t)}dt  $$
then, the resulting problem transforms to yet another formulation of a minimum-distance TSP (Cf.~Section 3). In fact, the most important aspect of \eqref{eq:prob-TSP-minTime-ocp-basic} is that it provides a clear avenue for generalization that is conducive to modeling traveling-salesman-type problems driven by constrained, nonlinear dynamical systems such as those illustrated in Fig.~\ref{fig:Grand-tour}.

\subsection{Generalization Based on Nonlinear Dynamics}
\begin{figure}[h!]
      \centering
      \framebox{\parbox{0.75\textwidth}{\centering
      {\includegraphics[trim={20mm, 25mm, 35mm, 15mm}, clip, angle=0, width = 0.65\textwidth]{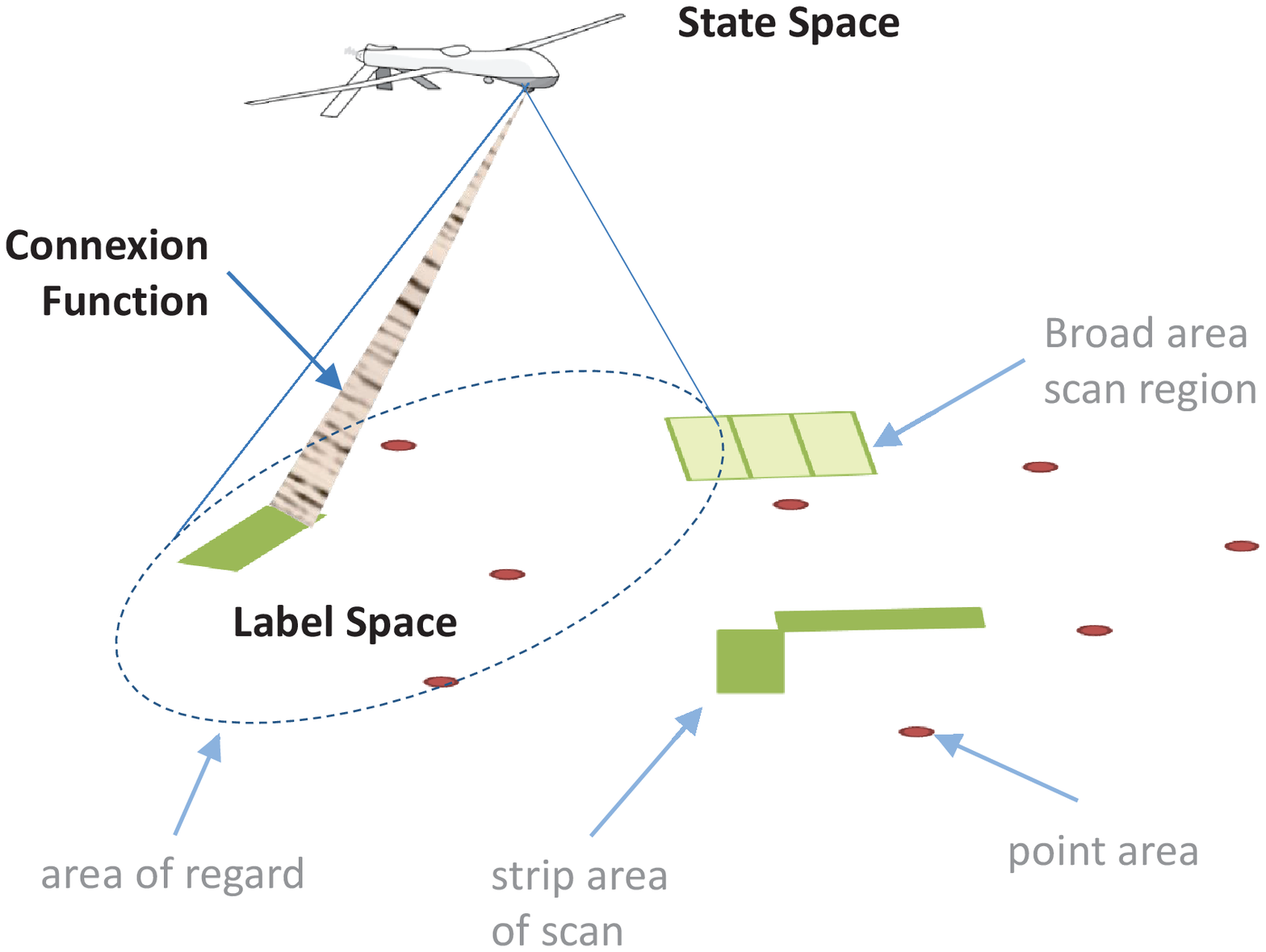}}
      }}
      \caption{\textsf{A conceptual illustration of state space, label space and the connexion function for a drone tasking problem\cite{TSP-UAV-ross}. }}
    \label{fig:UAV}
   \end{figure}

In many applications such as robotics \citep{TSP-Dubins}, the natural home space for the ``salesman'' is some state space $\X$ that is not necessarily the label space.  Furthermore, the dynamics of the salesman is given by some nonlinear controllable ordinary differential equation of the type,
\begin{equation}\label{eq:xdot=f}
\dot\bx = \bff(\bx, \bu, t)
\end{equation}
where, $\bx \in \real{N_x}$ is a state variable, $\bu \in \real{N_u}$ is the control variable and $\bff: \real{N_x} \times \real{N_u} \times \Real \to \real{N_x}$ is the nonlinear dynamics function. Because all aspects of the problem definition so far are described in terms of a label space, it is necessary to connect it to the state space. Let this connection be given by some algebraic equation,
\begin{equation}\label{eq:g=0-nou}
\bg(\bl, \bx, t) = \bzero
\end{equation}
where, $\bg: \real{N_l} \times \real{N_x} \times \Real \to \real{N_g} $ is called a connexion function. An example of a physical description of state space, label space and the connexion function is illustrated in Fig.~\ref{fig:UAV}.
%
%======================================================================================
%==========================================================================================
%
In this example, an uninhabited aerial vehicle (UAV) is tasked to collect over a geographic region \citep{TSP-UAV-shima,TSP-UAV-ross}.  The areas of interest vary from a point area to a broad area of scan. Technical properties associated with the areas of interest are defined in label space.  The state variables of the UAV is given in terms of of its position, velocity and the orientation of the maneuverable camera.  The connexion function is the mathematical model that connects the label space variables to the state space variables based on the precise position and orientation of the camera at a given instant of time.

From \eqref{eq:xdot=f} and \eqref{eq:g=0-nou}, it follows that \eqref{eq:ldot=w} may be written implicitly as a differential algebraic equation,
\begin{equation}\label{eq:ldot-via-xdot}
{\left.
                               \begin{array}{l}
                                 \dot\bx = \bff(\bx, \bu, t)  \\
                                 \bzero = \bg(\bl, \bx, t)
                               \end{array}
                             \right\} }  \quad \Longrightarrow \quad \dot\bl = \bw
\end{equation}
The significance of replacing \eqref{eq:ldot=w} by the state-space representation is that it facilitates a direct means to incorporate the full nonlinear dynamics of a salesman in the optimization problem formulation.

\subsection{Generalization Based on Economics-Driven Cost/Payoff/Constraint Functionals}
As a result of \eqref{eq:ldot-via-xdot}, the decision variables expand to the tuple, $(\bl(\cdot), \bx(\cdot), \bu(\cdot), t_0, t_f)$.  This can be taken to great advantage by formulating cost, payoff and/or constraint functionals in their more natural ``home spaces.''  For example, in the electric vehicle-routing problem \citep{TSP-Electric}, the capacity constraint of the battery is more naturally expressed as a functional constraint in terms of the state $\bx$ of the vehicle \citep{EV-battery},
\begin{equation}\label{eq:batt-state}
\int_{t_0}^{t_f}S(\bx(t))\,dt \ge C_{safe}
\end{equation}
In \eqref{eq:batt-state}, $S$ is the state-of-charge and $C_{safe}$ is a battery charge required for safe operations.  In the absence of extending the decision variables to incorporate the state vector, \eqref{eq:batt-state} would need to be transformed to a label space constraint using \eqref{eq:g=0-nou}.  This difficult task is completely circumvented by incorporating the state vector as part of the optimization variable. The alternative is to generate proxy models \citep{TSP-Electric} as a means to extend the scope of vehicle routing problems \citep{vidal}.

In \eqref{eq:prob-TSP-minTime}, the travel time is an implicit functional of the label-space trajectory $\bl(\cdot)$.  Replacing \eqref{eq:ldot=w} by \eqref{eq:ldot-via-xdot}, the travel time now becomes an implicit functional of $\bl(\cdot), \bx(\cdot)$ and $\bu(\cdot)$. That is, the functional,
\begin{equation}\label{eq:J=travel-time}
J_{time}[\bl(\cdot), \bx(\cdot), \bu(\cdot), t_0, t_f]:= t_f - t_0
\end{equation}
with constraints given by \eqref{eq:ldot-via-xdot} offers a more accurate model of travel time with obvious business implications in practical routing problems. Note that \eqref{eq:J=travel-time} also allows the clock time to be optimized.  Equation \eqref{eq:J=travel-time} may be further modified to take into account gas/power consumption depending upon the vehicle type (standard/hybrid/electric).  Analogous to \eqref{eq:batt-state}, gas/power consumption may be written in terms of an integral,
\begin{equation}\label{eq:burn-rate-integ}
\displaystyle \int_{t_0}^{t_f} f_0\big(\bx(t), \bu(t), t\big)\, dt
\end{equation}
where, $f_0$ is a function that models the time-varying gas/power consumption rate.  This function may be well-modeled using the physics of the automobile power train \citep{EV-battery}. A convex combination of \eqref{eq:burn-rate-integ} with \eqref{eq:J=travel-time} generates the functional,
\begin{multline}\label{eq:J=hybrid}
J_{hybrid}[\bl(\cdot), \bx(\cdot), \bu(\cdot), t_0, t_f]:= \alpha (t_f - t_0) \\
+ (1-\alpha) \displaystyle \int_{t_0}^{t_f} f_0\big(\bx(t), \bu(t), t\big)\, dt \qquad \alpha \in [0, 1]
\end{multline}
The parameter $\alpha$ in \eqref{eq:J=hybrid} offers a sliding scale over the ``trade-space'' of time and energy consumption.

\subsection{An Optimal Control Framework for a TSP and its Variants }
Substituting \eqref{eq:ldot-via-xdot} in \eqref{eq:prob-TSP-minTime-ocp-basic} while simultaneously adding additional levels of abstraction for the functionals associated with the vertices and arcs of the underlying $\mathcal{T}$-graph, we arrive at:
% The Traveling X-Problem
%
%==========================
\begin{eqnarray}
&\big(\mathcal{T_X}P\big) \left\{
\begin{array}{lll}
\displaystyle\mathop\text{Minimize}   & \displaystyle J[\bl(\cdot), \bx(\cdot), \bu(\cdot), t_0, t_f] \\
\text{Subject to }
& \dot\bx(t) = \bff(\bx(t), \bu(t), t)\\
& \bg(\bl(t), \bx(t), t) = \bzero \\
%&  \bu(t) \in \mathbb{U}(t, \bx(t))\\
& K^m\left[\bl(\cdot), \bx(\cdot), \bu(\cdot), t_0, t_f\right] \left\{
                              \begin{array}{ll}
                                \le 0, & \forall\ m \in N_{\ge 0} \\
                                = 0, & \forall\ m \in N_{=0}
                              \end{array}
                            \right. \\[0.1em]
& \big(\bl(t_0), \bl(t_f) \big) \in \Lspace^b \subseteq \Lspace \\
& \bx(t) \in \mathbb{X}(t)\\
& \bu(t) \in \mathbb{U}(t, \bx(t))
\end{array} \right.& \label{eq:prob-TXP}
\end{eqnarray}
%========================================
%
In \eqref{eq:prob-TXP}, the objective function is simply stated in terms of an abstract functional $J$ in order to facilitate the formulation of a generic cost functional that go beyond those discussed in Sec.~5.2.  Likewise, the functionals $K^m$ in \eqref{eq:prob-TXP} may be time-on-task functionals, walk-indicator functionals, degree functionals, capacity-constraint functionals or some other functionals. Furthermore, because each vertex may have several constraints, the total number of these constraints may be greater than $\abs{N_v}$. The index sets $N_{=0}$ and $N_{\ge 0}$ simply organize the constraint functionals into equalities and inequalities.  The set $\Lspace^b \subseteq \Lspace $ stipulates the boundary conditions for the label space trajectory. Also included in \eqref{eq:prob-TXP} are state variable constraints $\bx(t) \in \X(t)$ and state-dependent control constraints given by $\bu(t) \in \mathbb{U}(t, \bx(t))$.  Such constraints are included in Problem $(\mathcal{T_X}P) $  because they are critically important in practical applications\cite{ross-book,vinter} as well as coordinate transformations of optimal control problems\cite{ross-book}.

Because \eqref{eq:prob-TXP} is a generalization of the problems discussed in Sections 3--5, it is evident that a fairly large class of traveling-salesman-like problems can be modeled as instances of Problem~$(\mathcal{T_X}P)$.

\section{Illustrative Numerical Examples}
A motorized TSP was first introduced in \cite{TSP-motorized}. In essence, a motor is a differential equation; hence, a motorized TSP is one with dynamical constraints.  Here, we combine this problem with several other variants of the TSP \citep{CETSP,TSPFN} to generate a minimum-time close-enough motorized traveling salesman problem with forbidden neighborhoods given by:
%
%==========================
\begin{eqnarray}
& \bl \in \real{4}, \quad \bx \in \real{4}, \quad \bu \in \real{2} \nonumber\\
&\big(\textit{fastCEMTSPFN-1}\big) \left\{
\begin{array}{lll}
\displaystyle\mathop\text{Minimize}   & \displaystyle J[\bl(\cdot), \bx(\cdot), \bu(\cdot), t_0, t_f]:= t_f-t_0 \\ %[0.5em]
\text{Subject to }
& \dot x_1(t) = x_3(t), \quad \abs{x_3(t)} \le 1 \quad \forall\ t \in [t_0, t_f] \\
& \dot x_2(t) = x_4(t), \quad \abs{x_4(t)} \le 1 \quad \forall\ t \in [t_0, t_f] \\
& \dot x_3(t) = u_1(t), \quad \abs{u_1(t)} \le 1 \quad \forall\ t \in [t_0, t_f] \\
& \dot x_4(t) = u_2(t), \quad \abs{u_2(t)} \le 1 \quad \forall\ t \in [t_0, t_f] \\
& l_i(t)-x_i(t) = 0,  \quad \forall\ t \in [t_0, t_f]\\
& \qquad i = 1, 2, 3, 4 \\
& T^i\left[\bl(\cdot), t_0, t_f\right] \ge \tau^i_{min}, \quad \forall\ i = 1, \ldots, N_v \\[0.5em]
& \displaystyle \left( \frac{x_1(t) - x_1^i}{a^i} \right)^2 + \left( \frac{x_2(t) - x_2^i}{b^i} \right)^2 \ge 1,\\
& \qquad \forall\ i = 1, \ldots N_{obs} \\
& \big(\bl(t_0), \bl(t_f) \big) = (\bzero, \bzero)
\end{array} \right.& \label{eq:prob-EX-1}
\end{eqnarray}
%========================================
%
It is apparent that \eqref{eq:prob-EX-1} is a special case of \eqref{eq:prob-TXP}. In fact, it is a generalization of the fast TSP posed in \eqref{eq:prob-TSP-minTime-ocp-basic}. Equation \eqref{eq:prob-EX-1} is motorized by the four differential equations resulting from an elementary application of Newtonian dynamics: The position and velocity of the salesman are $(x_1, x_2) \in \real{2}$ and $(x_3, x_4) \in \real{2}$ respectively. Both the velocity and acceleration $(u_1, u_2) \in \real{2}$ vectors are constrained in the $\ell_\infty$ norm (i.e., components are constrained in terms of their absolute values).  The $N_{obs}$ ellipsoidal constraints in \eqref{eq:prob-EX-1} are the forbidden neighborhoods whose centroids are given by $(x_1^i, x_2^i), \ i = 1, \ldots, N_{obs}$.  The parameters $a^i > 0$ and $b^i > 0$ determine the shape of the ellipse. The last constraint in \eqref{eq:prob-EX-1} requires the salesman to start and end at the origin.

\begin{figure}[h!]
      \centering
      \framebox{\parbox{0.75\textwidth}{\centering
      {\includegraphics[trim={1mm, 1mm, 1mm, 1mm}, clip, angle=0, width = 0.63\textwidth]{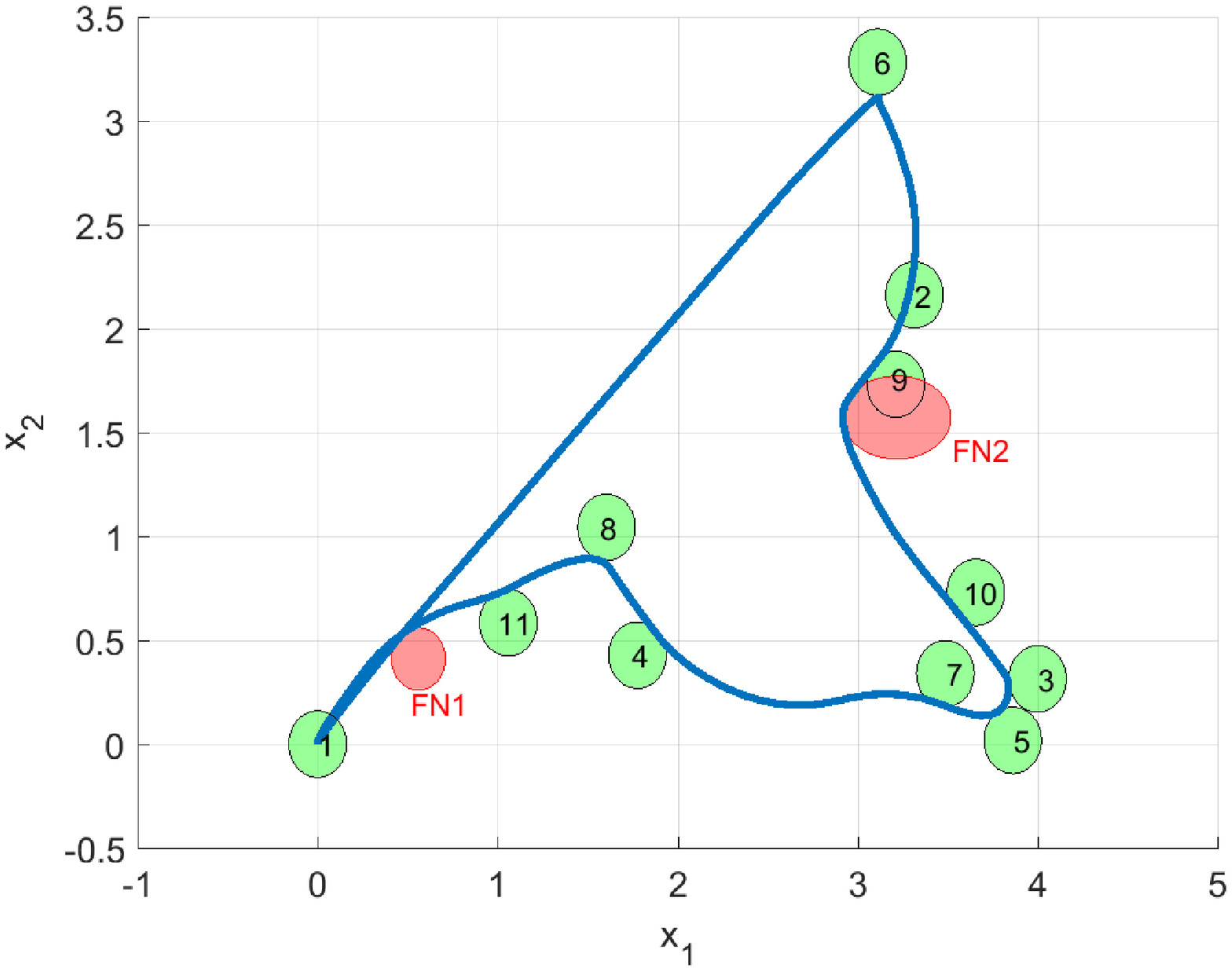}}
      }}
      \caption{\textsf{Solution to Problem (fastCEMTSPFN-1) for 11 cities and 2 forbidden neighborhoods marked FN1 and FN2.}}
    \label{fig:TSPResult11BoxFN}
   \end{figure}

A sample data set for \eqref{eq:prob-EX-1} with $N_v = 11$ and $N_{obs} = 2$ is shown in Fig.~\ref{fig:TSPResult11BoxFN}.
%
%======================================================================================
%==========================================================================================
%
Also shown in Fig.~\ref{fig:TSPResult11BoxFN} is the solution obtained by solving \eqref{eq:prob-EX-1} using DIDO's implementation\cite{ross-dido} of pseudospectral optimal control theory \citep{PSOC}. Some key points to note regarding the solution presented in Fig.~\ref{fig:TSPResult11BoxFN} are:
\begin{enumerate}
\item Nearly all the ``arcs'' between city-pairs are curvilinear. This is because the salesman's trajectory is required to satisfy the Newtonian dynamical constraints as well as the instantaneous $\ell_\infty$ bounds on velocity and acceleration as dictated in \eqref{eq:prob-EX-1}.
\item The city-pair arcs as well as the entry and exit points to the city are a natural outcome of solving the posed problem (Cf.~\eqref{eq:prob-EX-1}).  That is, the city-pair arcs were not determined either \textit{a priori} or \textit{a posteriori} to the determination of the city sequence.  See also Remark \ref{rem:newTSP-objFun}.
\item The ellipsoidal forbidden neighborhood marked FN2 overlaps neighborhood No.~9. Thus, the allowable region for neighborhood No.~9 is nonconvex.
\end{enumerate}

\begin{figure}[h!]
      \centering
      \framebox{\parbox{0.75\textwidth}{\centering
      {\includegraphics[trim={1mm, 1mm, 1mm, 1mm}, clip, angle=0, width = 0.65\textwidth]{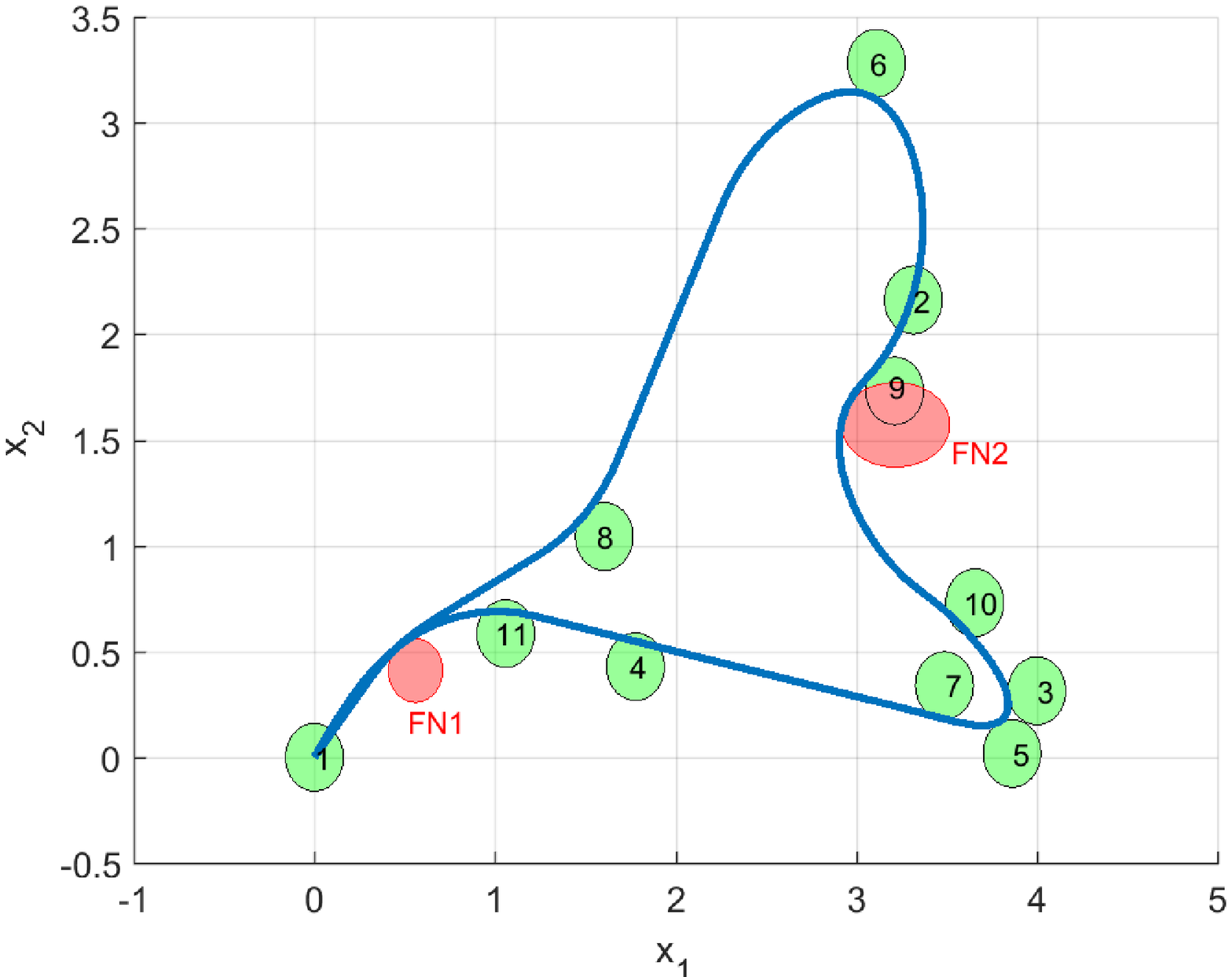}}
      }}
      \caption{\textsf{Solution to Problem (fastCEMTSPFN-2) for the same data shown in Fig.~\ref{fig:TSPResult11BoxFN}. }}
    \label{fig:TSPResult11DiskFN}
   \end{figure}

Next, consider a ``variant'' of \eqref{eq:prob-EX-1} obtained by replacing the $\ell_\infty$ constraints on velocity and acceleration by its $\ell_2$ version,
\begin{equation}\label{eq:prob-EX-2}
(\textit{fastCEMTSPFN-2}): \quad \sqrt{x_3^2(t) + x_4^2(t)} \le 1, \ \sqrt{u_1^2(t) + u_2^2(t)} \le 1 \qquad \forall\ t \in [t_0, t_f]
\end{equation}
That is, Problem \textit{(fastCEMTSPFN-2)} is identical to Problem \textit{(fastCEMTSPFN-1)} except for the additional nonlinear constraint given by \eqref{eq:prob-EX-2}.  A solution to this problem is shown in Fig.~\ref{fig:TSPResult11DiskFN}.
%
%======================================================================================
%==========================================================================================
%
It is apparent that Fig.~\ref{fig:TSPResult11DiskFN} differs from Fig.~\ref{fig:TSPResult11BoxFN} in several ways:
\begin{enumerate}
\item The turn at neighborhood No.~6 in Fig.~\ref{fig:TSPResult11BoxFN} is significantly sharper than the corresponding one in Fig.~\ref{fig:TSPResult11DiskFN}. The curvilinear turn  in Fig.~\ref{fig:TSPResult11DiskFN} is a clear and direct result of the $\ell_2$ constraint given by \eqref{eq:prob-EX-2}.
\item The sequence of visits in Fig.~\ref{fig:TSPResult11DiskFN} is different from that of Fig.~\ref{fig:TSPResult11BoxFN}; for example, compare the visit to neighborhood No.~8.
\item The entry and exit points to the same numbered cities in Figs.~\ref{fig:TSPResult11BoxFN} and \ref{fig:TSPResult11DiskFN} are different; for example, compare neighborhoods Nos.~8, 11 and 4.
\end{enumerate}

\begin{figure}[h!]
      \centering
      \framebox{\parbox{0.75\textwidth}{\centering
      {\includegraphics[trim={1mm, 1mm, 1mm, 1mm}, clip, angle=0, width = 0.65\textwidth]{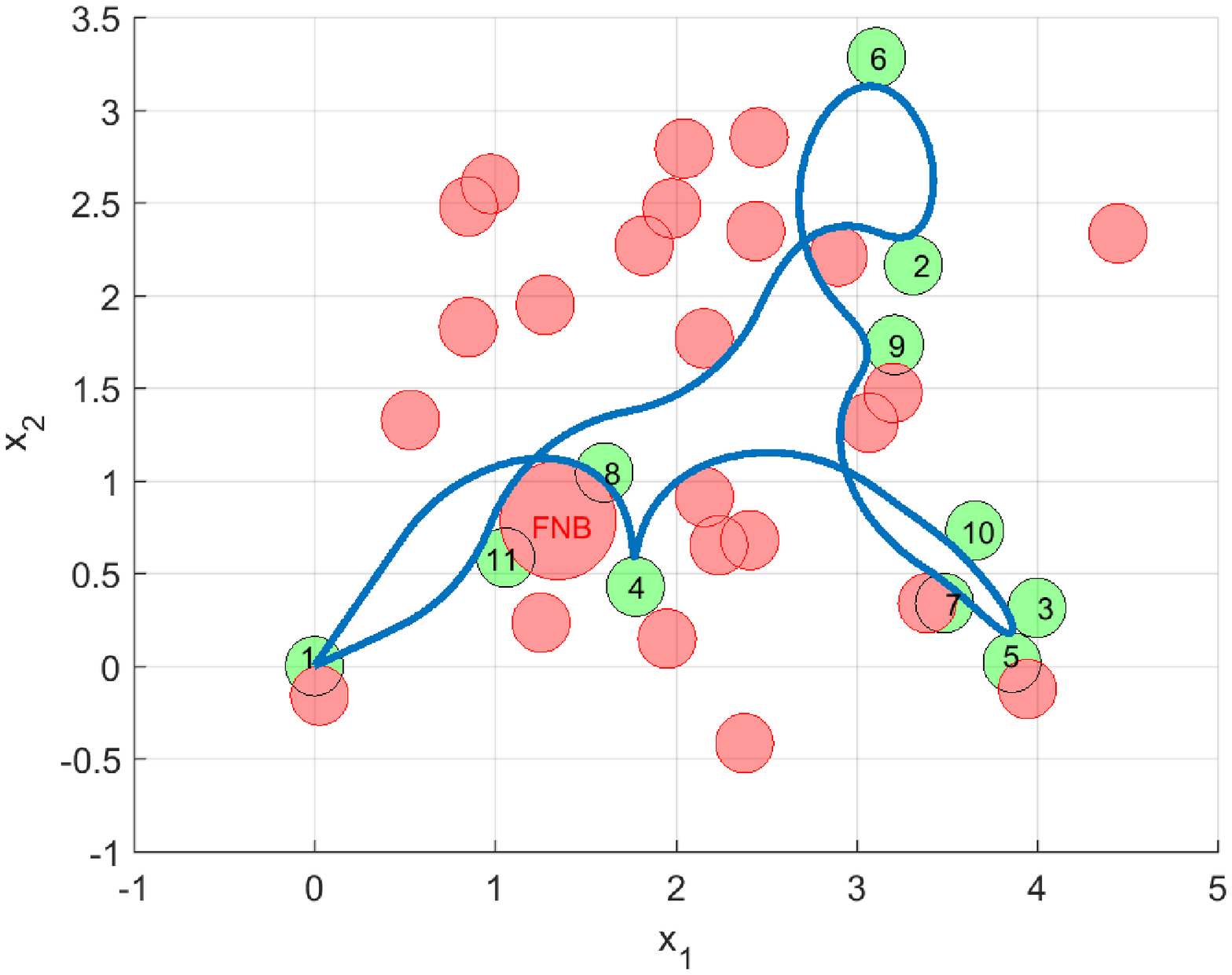}}
      }}
      \caption{\textsf{Solution to Problem (fastCEMTSPFN-2) for a new data set of nonconvex allowable and forbidden neighborhoods. }}
    \label{fig:TSPFNgaloreNo1}
   \end{figure}

For the third and final case, we change the problem ``data set'' with a random distribution of a large number of forbidden zones as shown in Fig.~\ref{fig:TSPFNgaloreNo1}. Only the neighborhood marked ``FNB'' was not randomly selected; rather, its size and location were purposefully set as indicated in Fig.~\ref{fig:TSPFNgaloreNo1}  to generate a more interesting case study.
%
%======================================================================================
%==========================================================================================
%
Also shown in Fig.~\ref{fig:TSPFNgaloreNo1} is a solution to the problem whose vehicle is constrained by \eqref{eq:prob-EX-2}.  Beyond the apparent efficacy of the process in obtaining a solution, noteworthy points pertaining to Fig.~\ref{fig:TSPFNgaloreNo1} are:
\begin{enumerate}
\item Some of the forbidden neighborhoods are nonconvex.  The nonconvex sets are generated randomly in the sense that the circular neighborhoods are allowed to overlap.
\item Some of the allowable neighborhoods are nonconvex because the forbidden circular disks overlap the allowable regions.
\item Comparing Figs.~\ref{fig:TSPFNgaloreNo1} and \ref{fig:TSPResult11DiskFN}, it is clear that the addition of forbidden neighborhoods has drastically altered the sequence of visits.
\end{enumerate}

\section{Conclusions}

Many types of objective functions and constraints in emerging variations of the traveling salesman problem (TSP) are more naturally defined in their continuous-time home space.  Modeling these variants of the TSP using a classic discrete optimization framework is neither straightforward nor easy.  By inverting the traditional modeling process, that is, by formulating the naturally discrete quantities in terms of continuous-time nonsmooth functions, it is possible to generate a new framework for the TSP and some of its variants.  In this context, the discrete-optimization formalism of a TSP may be viewed as a problem formulation in the co-domain of the functionals that constitute its $\mathcal{T}$-graph.  In sharp contrast, the problem formulation presented in this paper is in the domain of the functionals of the TSP $\mathcal{T}$-graph.
%In terms of the notion of an $\mathcal{F}$-graph, the problem formulations presented in this paper are domain-centric while the discrete-optimization formalism is co-domain-centric.
This perspective suggests that the discrete-optimization- and the variational formulations are effectively two sides of the same $\mathcal{T}$-graph constructs introduced in this paper.

There is no doubt that there are a vast number of open theoretical issues in the proposed framework.  Nonetheless, we have demonstrated, by way of a fast, close-enough motorized TSP with forbidden neighborhoods, that it is indeed possible to solve some challenging problems using the new formalisms. One of the interesting revelations indicated by the numerical studies is the big impact of seemingly small changes in the motion-constraints of the traveling salesman. This suggests that an exploitation of nonlinear models may provide a discriminating edge to a business/engineering operation by way of a non-intuitive economical utilization of its end-to-end systems.

\section*{Acknowledgments}
Funding from the Defense Advanced Research Projects Agency, the Center for Multi-INT Studies, and the Department of Defense (DoD) are gratefully acknowledged. The views, opinions and/or findings expressed are those of the authors and should not be interpreted as representing the official views or policies of the DoD or the U.S. Government.

Strong support from the Secretary of the Navy and United States Navy's Judge Advocate General's Office is gratefully acknowledged. Portions of this work are protected by the following patents issued by the United States Patent and Trademark Office: Patent Nos. US 9,983,585 B1; US 10,476,584 B1; and 62/971,068 (patent pending).
%\bibliography{mybibfile}

\begin{thebibliography}{10}


\bibitem[{Aggarwal et~al.(1999)}]{TSP-angle}
Aggarwal A, Coppersmith D, Khanna, S, Motwani R, Schieber B (1999) The angular-metric traveling salesman problem.
{\it SIAM J. Comput.} 29(3):697--711.


\bibitem[{Arkin and Hassin(1994)}]{TSPN}
Arkin EM, Hassin R (1994) Approximation algorithms for the goemetric covering salesman problem.
{\it Discrete Appl. Math.} 55(3):197--218.


\bibitem[{Chowdhury et~al.(2019)}]{DTSP-2019}
Chowdhury S, Marufuzzaman M, Tunc H, Bian L, Bullington W (2019) A modified ant colony optimization algorithm to solve a dynamic traveling salesman problem: A case study with drones for wildlife surveillance. {\it J. Comp. Design and Engr.} 6(3):368--386.

\bibitem[{Cook(2012)}]{cook-2012}
Cook WJ (2012) {\it In Pursuit of the Traveling Salesman: Mathematics at the Limits of Computation} (Princeton University Press, Princeton, NJ).

\bibitem[{Dantzig, Fulkerson and Johnson(1954)}]{Dantzig-54}
Dantzig G, Fulkerson R, Johnson S (1954) Solution of a large-scale traveling-salesman problem. {\it J. Oper. Res. Soc. of America} 2(4):393--410.

\bibitem[{Fischer and Hungerlander(2017)}]{TSPFN}
Fischer A, Hungerlander P (2017) The traveling salesman problem on grids with forbidden neighborhoods.
{\it J. Combinatorial Optimization} 34(3):891--915.

\bibitem[{Flood(1956)}]{Flood}
Flood MM (1956) The traveling-salesman problem. {\it Oper. Res.} 4(1):61--75.

\bibitem[{Gentilini, Margot and Shimada(2013)}]{TSPN-robotics}
Gentilini I, Margot F, Shimada K (2013) The travelling salesman problem with neighbourhoods: MINLP solution. {\it Optimization Methods and Software} 28(2): 364--378.

\bibitem[{Gottlieb and Shima(2015)}]{TSP-UAV-shima}
Gottlieb Y, Shima T (2015) UAVs task and motion planning in the presence of obstacles and prioritized targets. {\it Sensors} 15:29734--29764.

\bibitem[{Gulczynski, Heath and Price(2006)}]{CETSP}
Gulczynski DJ, Heath JW, Price CC (2006) The close enough traveling salesman problem: A discusion of several heuristics. Alt FB, Fu MC Golden BL ed. {\it Perspectives in Op. Res.: Papers in Honor of Saul Gass' 80th Birthday } (Springer, Boston)  271--283

\bibitem[{Gunawan, Hoong and Vansteenwegen(2016)}]{OP-survey}
Gunawan A, Hoong CL, Vansteenwegen P (2016) Orienteering Problem: A survey of recent variants, solution approaches and applications. {\it European J. Op. Res.} 255:315--332.

\bibitem[{Gutin and Punnen(2007)}]{TSP-variants}
Gutin G, Punnen AP (2007) {\it The Traveling Salesman Problem and its Variations} (Springer, New York)

\bibitem[{Ma and Castanon(2006)}]{TSP-Dubins}
Ma X, Castanon DA (2006) Receding horizon planning for Dubins traveling salesman problems. {\it Proc. 45th IEEE CDC} 5453--5458.

\bibitem[{Psaraftis(1988)}]{DTSP}
Psaraftis, HN (1988) Dynamic vehicle routing problems. B. Golden, A. Assad (Editors), Vehicle routing: methods and studies, (Elsevier Science Publishers BV) 223--248.

\bibitem[{Miller, Tucker and Zemlin(1960)}]{MTZ}
Miller CE, Tucker AW, Zemlin RA (1960) Integer programming formulations and traveling salesman problems. {\it J. ACM} 7, 326--329.

\bibitem[{Nesterov(1983)}]{nesterov83}
Nesterov YE A method of solving a convex programming problem with convergence rate $\mathcal{O}(1/k^2)$, \textit{Soviet Math. Dokl.}, 27/2 371--376 (Translated by A. Rosa).

\bibitem[{Polyak(1964)}]{polyak64}
Polyak BT (1964) Some methods of speeding up the convergence of iteration methods, \textit{USSR Computational Math. and Math. Phys.}, 4/5, 1--17 (Translated by H. F. Cleaves).

\bibitem[{Rana, Anand and Bose(2019)}]{TSP-telescope}
Rana J, Anand S, Bose S (2019) Optimal search strategy for finding transients in large-sky error regions under realistic constraints. {\it The Astrophysical Journal } 876(2):104.

\bibitem[{Ross(2020)}]{ross-dido}
Ross IM (2020) Enhancements to the DIDO$^\copyright$ optimal control toolbox. arXiv preprint.  arXiv:2004.13112.

\bibitem{ross-accel}
Ross IM (2019) An optimal control theory for accelerated optimization. arXiv preprint arXiv:1902.09004.

\bibitem[{Ross(2019)}]{ross-jcam-1}
Ross IM (2019) An optimal control theory for nonlinear optimization. {\it J. Computational and Appl. Math.} 354:39--51.


\bibitem{ross-book}
Ross IM (2015) \textit{A Primer on Pontryagin's Principle in Optimal Control}, Second Edition, Collegiate Publishers, San Francisco, CA.

\bibitem[Ross and Karpenko(2012)]{PSOC}
Ross IM, Karpenko M (2012) A review of pseudospectral optimal control: From theory to flight. {\it Annual Reviews in Control} 36:182--197.

\bibitem[{Ross, Proulx and Karpenko(2019)}]{TSP-UAV-ross}
Ross IM, Proulx RJ, Karpenko M (2019) Autonomous UAV sensor planning, scheduling and maneuvering: An obstacle engagement technique. {\it ACC} 65--70.


\bibitem[{Schneider, Stenger and Goeke(2014)}]{TSP-Electric}
Schneider M, Stenger A, Goeke D (2014) The electric vehicle-routing problem with time windows and recharging stations. {\it Transportation Sc.} 48(4): 500--520.

\bibitem[{Sciarretta, Back and Guzzella(2004)}]{EV-battery}
Sciarretta A, Back M, Guzzella L (2004) Optimal control of parallel hybrid electric vehicles. {\it IEEE Tran. Control Sys. Tech.} 12(3): 352--363.

\bibitem[von Stryk and Glocker(2001)]{TSP-motorized}
von Stryk O, Glocker M (2001) Numerical mixed-integer optimal control and motorized traveling salesmen problems. {\it J. Europ\'{e}en des Syst\`{e}mes Automatis\'{e}s} 35(4): 519--533.

\bibitem[{Vinter(2000)}]{vinter}
Vinter RB (2000) {\it Optimal Control} (Birkh\"{a}user, Boston)

\bibitem[{Vidal(2017)}]{vidal}
Vidal T (2017) Node, edge, arc routing and turn penalties: Multiple problems--one neighborhood extension. {\it Op. Res.} 65(4): 992--1010.

\bibitem[{Wang, Golden and Wasil(2019)}]{wang2019}
Wang X, Golden B, Wasil E (2019) A Steiner zone variable neighborhood search heuristic for the close-enough traveling salesman problem. {\it Computers and Op. Res.} 101: 200--219.

\bibitem[{Witze(2018)}]{witze}
Witze A (2018) Jupiter has 10 more moons we didn't know about --- and they're weird. {\it Nature} 559:312--313.


\end{thebibliography}

\end{document}